\documentclass{birkjour}
 \theoremstyle{definition}
 
 \theoremstyle{remark}

 \numberwithin{equation}{section}

\usepackage{array}

\makeatletter
\newcommand{\thickhline}{%
    \noalign {\ifnum 0=`}\fi \hrule height 1pt
    \futurelet \reserved@a \@xhline
}
\newcolumntype{"}{@{\hskip\tabcolsep\vrule width 1pt\hskip\tabcolsep}}
\makeatother

\usepackage{amsmath,amsthm,amsxtra}

\usepackage{MnSymbol}

\usepackage{stmaryrd}

\usepackage{graphics}
\usepackage{psfrag}
\usepackage{enumerate}
\usepackage[format=hang,indention=0pt,figurename=Fig.]{caption}
 \usepackage[dvips]{color}

\definecolor{refkey}{gray}{0.5}
\definecolor{labelkey}{gray}{0.5}

\newtheorem{theorem}{Theorem}[section]
\newtheorem{proposition}[theorem]{Proposition}
\newtheorem{lemma}[theorem]{Lemma}
\newtheorem{corollary}[theorem]{Corollary}

\theoremstyle{definition}

\newtheorem{example}[theorem]{Example}

\newtheorem{remark}[theorem]{Remark}

\numberwithin{equation}{section}

\newcommand{\cQ}{\mathcal Q}

\newcommand{\CT}{\mathcal{C}}

\newcommand{\nR}{\mathbb{R}}
\newcommand{\nM}{\mathbb{M}}
\newcommand{\nN}{\mathbb{N}}
\newcommand{\nI}{\mathbb{I}}
\newcommand{\nS}{\mathbb{S}}
\newcommand{\nT}{\mathbb{T}}
\newcommand{\nRc}{{\mathbb R}\cup\{\infty\}}

\newcommand{\lep}{\alpha}
\newcommand{\rep}{\beta}
\newcommand{\fvk}{\kappa}

 \DeclareMathSymbol{\varkappa} {\mathord}{AMSb}{"7B}

\newcommand{\Bo}{\medsquare}
\newcommand{\BoI}{\meddiamond}
\newcommand{\Ton}[1]{#1^{\boxminus}}
\newcommand{\TonI}[1]{#1^{\diamondminus}}
\newcommand{\pn}[1]{\llparenthesis #1 \rrparenthesis}

\newcommand{\Vs}[1]{\mathsf{#1}}
\newcommand{\Vb}[1]{\mathbf{{#1}}}

 \newcommand{\pI}{\mathsf{p}}

 \DeclareMathOperator{\Rot}{Rot}
 \DeclareMathOperator{\Ref}{Ref}

\begin{document}

\title[Triangles and groups via cevians]
 {Triangles and groups via cevians}

\author[\'A. B\'enyi]{\'Arp\'ad B\'enyi}

\address{%
Department of Mathematics \\
Western Washington University \\
516 High Street, Bellingham \\
Washington 98225, USA}

\email{Arpad.Benyi@wwu.edu}

\thanks{This work is partially supported by a grant from the Simons Foundation (No.~246024 to \'Arp\'ad B\'enyi).}

\author[B. \'{C}urgus]{Branko \'{C}urgus}
\address{%
Department of Mathematics \\
Western Washington University \\
516 High Street, Bellingham \\
Washington 98225, USA}

\email{Branko.Curgus@wwu.edu}



\date{\today}


\begin{abstract}
For a given triangle $T$ and a real number $\rho$ we define Ceva's triangle $\CT_\rho(T)$ to be the triangle formed by three cevians each joining a vertex of $T$ to the point which divides the opposite side in the ratio $\rho:(1-\rho)$.  We identify the smallest interval $\nM_T \subset \nR$ such that the family $\CT_\rho(T), \rho\in \nM_T$, contains  all Ceva's triangles up to similarity. We prove that the composition of operators $\CT_\rho, \rho \in \nR$, acting on triangles is governed by a certain group structure on $\nR$. We use this structure to prove that two triangles have the same Brocard angle if and only if a congruent copy of one of them can be recovered by sufficiently many iterations of two  operators $\CT_\rho$ and $\CT_\xi$ acting on the other triangle.
\end{abstract}

 \maketitle


\section{Introduction}

A {\em median} of a triangle is a line segment connecting a vertex
to the midpoint of the opposite side. To each vertex of a triangle
corresponds exactly one median. A classical theorem in triangle
geometry states that the three medians of a given triangle form a
triangle. This new triangle is called the {\em median triangle}.
Moreover, the area of the median triangle is 3/4 of the area of the
host triangle. Two existence proofs of the median triangle and a
connection to a Heron-type formula for medians are revisited in
\cite{Benyi}; see also \cite{Hajja1} and the references therein.
A less well known result is that the median triangle of the median triangle is similar to the given triangle in the ratio $3/4$. This property is recalled by Scott in \cite{Scott} where it is referred to as the binary similarity of the sequence of median triangles.

The binary similarity property of the sequence of median triangles
was reformulated and extended by Griffiths \cite{Griffiths} as a
statement about a special class of linear operators mapping the
three dimensional Euclidean space into itself. Griffiths refers to
these operators as being of cyclically symmetrical type; they are
simply those operators whose matrix representations with respect to
the standard basis are left-circulant matrices with orthogonal rows. This nice connection with matrix algebra can be used to produce an infinite number of such binary sequences, see \cite[Proposition~2]{Griffiths}. We will return to this observation shortly, since it turns out to be closely connected to our work.

A median of a triangle is just a special cevian. For a given triangle $T$ and a real number $\rho$, instead of the medians we can consider three cevians, each joining a vertex of $T$ to the point which divides the opposite side in the ratio $\rho:(1-\rho)$. Three such cevians also form a triangle. With a particular choice of order of these cevians, we call the triangle formed in  this way Ceva's triangle of $T$ and denote it by $\CT_\rho(T)$. With a different order of sides, such triangles are considered by Hajja in \cite{Hajja1, Hajja2}, where they are called  $s$-median or generalized median triangles. By analogy with the iterative procedure for median triangles that leads to binary similarity, it is natural to ask whether the same holds for other sequences of nested triangles. 
The limiting behavior, in shape, of various nested sequences of triangles
constructed iteratively was considered by many authors; more recently by
Ismailescu and Jacobs in \cite{Ismailescu} where one can find other references. This work has motivated Hajja \cite{Hajja1, Hajja2} to do the same for the sequence of generalized median triangles. Using a suitable shape function written in terms of the side lengths of the original triangle, \cite[Theorem~3.1]{Hajja2} reveals a delicate limiting behavior of the sequence of generalized median triangles. Similar results were obtained for a related iteration process on triangles in the remarkable paper of Nakamura and Oguiso \cite{Nakamura} by using moduli space of the similarity classes of triangles. And so, it might seem that this is the end of the story as far as the iterated sequence of generalized median triangles is concerned. However, returning to Griffiths' observation regarding cyclical symmetry, it turns out that the intricate behavior of the sequence of generalized median triangles introduced in \cite{Hajja1, Hajja2} depends precisely on one detail of its definition. Indeed, the order of the sides matters in the iteration process.  With our definition, we do have the binary similarity property of the sequence of Ceva's triangles, which is  completely  analogous to that of the median triangles.

The goal of this article is to provide a complete understanding of
the family of Ceva's triangles. A common thread throughout this work is the presence of a special group structure on the extended real line related to the family of Ceva's triangles. We elaborate on the observation made in \cite{Griffiths} and explain how linear algebra is connected to the group structure we alluded to before. These connections allow us, for example, to show that we can iterate Ceva's triangles with different parameters $\rho$ and calculate the parameter of the new Ceva's triangle obtained this way. We prove that the family $\CT_\rho(T), \rho\in [0,1)$, contains all Ceva's triangles of $T$ up to direct similarity. In fact, we can identify the smallest interval $\nM_T \subset [0,1)$ such that the family $\CT_\rho(T), \rho\in \nM_T$, contains  all Ceva's triangles up to similarity.  Incidentally, we also discover a new shape function which is closely related to the one introduced by Hajja in \cite{Hajja2}, as well as a characterization of the equality of the Brocard angles of two triangles in terms of their respective families of Ceva's triangles which extends the one in \cite[Theorem 3.1]{Hajja2}. This characterization of equality of the Brocard angles of two triangles is closely related to a theorem of Stroeker \cite[page~183]{Stroeker}. Lastly, returning to the generic question about the behavior of some iterative geometric process, we prove that, given two triangles having the same Brocard angle, we can recover a congruent copy of one of them by a sufficiently long iteration of two Ceva's operators acting on the other triangle.

\section{Basic notions} \label{sbn}

A {\em triangle} is a set of three noncollinear points and the three
line segments joining each pair of these points.  The triangle
determined by three noncollinear points $A, B, C$ is denoted by
$ABC$. The points $A, B, C$ are called {\em vertices} and the line
segments $a = BC, b = C\!A, c = AB$ are called {\em sides} of the
triangle $ABC$. The notation used for line segments will also stand for their   lengths. In particular, symbols $a, b, c$ denote the lengths
of the corresponding sides as well.  We will always label the vertices of a triangle {\em counterclockwise}. This convention is essential in the definition of Ceva's triangle in the next section. Also, in this way the sides become the oriented line segments $\overrightarrow{BC}, \overrightarrow{C\!A}, \overrightarrow{AB}$.  A triangle with such imposed orientation we call an {\em oriented triangle}. Notice that the counterclockwise order is also imposed on the lengths of the sides which we write as an ordered triple $(a,b,c)$ of positive real numbers. The adjective oriented will be omitted if it is clearly implied by the context in which a related triangle appears.

An ordered triple $(u,v,w)$ is {\em increasing} ({\em decreasing})
if $u < v < w$ (respectively, $u > v > w$). If $ABC$ is a scalene triangle with side lengths $a, b, c$, then we have the following dichotomy: the set
\begin{equation} \label{eqinde}
\bigl\{ (a,b,c), (b,c,a), (c,a,b) \bigr\}
\end{equation}
either contains a decreasing or an increasing triple. To justify
this, we can assume that $a=\min\{a,b,c\}$. Then, $c > b$ or  $b >
c$. If $c > b$, then $(a,b,c)$ is increasing. If  $b > c$, then
$(b,c,a)$ is decreasing. A scalene oriented triangle $ABC$ for which the set
\eqref{eqinde} contains an increasing (decreasing, respectively)
triple is called an {\em increasing} ({\em decreasing}) {\em
triangle}. For non-equilateral isosceles triangles we introduce the following intuitive terminology: if its legs are longer than its base we call it a {\em narrow} triangle; if its legs are shorter than its base we call it a {\em wide} triangle. For two oriented non-equilateral triangles we say that they have the {\em same orientation} if they are both increasing, or they are both decreasing, or they are both wide, or they are both narrow.

We recall the definitions of similarity and congruence for oriented triangles. Oriented triangles $ABC$ and $XY\!Z$ are {\em directly similar} if
\[
\frac{a}{x} = \frac{b}{y} = \frac{c}{z} \quad \text{or} \quad
 \frac{a}{y} = \frac{b}{z} = \frac{c}{x} \quad \text{or} \quad
 \frac{a}{z} = \frac{b}{x} = \frac{c}{y}.
\]
If $ABC$ and $XY\!Z$ are directly similar with $a/x = b/y = c/z = l$, then we will write $(a,b,c) = l (x,y,z)$.

Oriented triangles $ABC$ and $XY\!Z$ are {\em reversely similar} if
\[
\frac{a}{z} = \frac{b}{y} = \frac{c}{x} \quad \text{or} \quad
 \frac{a}{y} = \frac{b}{x} = \frac{c}{z} \quad \text{or} \quad
 \frac{a}{x} = \frac{b}{z} = \frac{c}{y}.
\]
Two oriented triangles are said to be {\em similar} if they are either directly or reversely similar. The common ratio of sides of two similar triangles is called the {\em ratio of similarity}. If the ratio of similarity is $1$, then directly (reversely, respectively) similar triangles are said to be {\em directly} ({\em reversely}) {\em congruent}. Notice that a triangle and its reflection are reversely congruent.

\section{Ceva's triangles} \label{srct}

Let $ABC$ be an oriented triangle and let $\rho$ be a real number.
Define the points $A_\rho, B_\rho$ and $C_\rho$ on the lines $BC,
C\!A$, and $AB$, respectively, by
\[
\overrightarrow{AC}_\rho = \rho \, \overrightarrow{AB}, \quad
\overrightarrow{C\!B}_\rho = \rho \, \overrightarrow{C\!A} , \quad
 \text{and} \quad
 \overrightarrow{BA}_\rho = \rho \, \overrightarrow{BC}.
\]
When $\rho\in (0, 1)$, the point $A_\rho$ is in the interior of the
line segment $BC$ while the cases $\rho > 1$ and $\rho < 0$ refer to
positions of the point exterior to the line segment $BC$. Also, $C_0
= A, B_0 = C, A_0 = B$ and $C_1 = B, B_1 = A, A_1 = C$. A similar
comment applies to the points $B_\rho$ and $C_\rho$. In this way we
obtain three {\it cevians}: $AA_\rho$, $BB_\rho,$ and $CC_\rho$. For
$\rho = 1/2$, they are medians.

For an oriented triangle $ABC$ and for an arbitrary $\rho \in {\mathbb R}$, the cevians $CC_\rho$, $BB_\rho$, and $AA_\rho$ form a triangle, see
\cite[Theorem~3.3]{Hajja1} and \cite[Theorem~2.7]{Hajja2}. Here is a
different, simple proof using vector algebra. Define the vectors
$\Vb{a} = \overrightarrow{BC}$, $\Vb{b}  = \overrightarrow{C\!A}$, and
$\Vb{c}  = \overrightarrow{AB}$ and $\Vb{x}_\rho =
\overrightarrow{CC_\rho}$, $\Vb{y}_\rho = \overrightarrow{BB_\rho}$,
$\Vb{z}_\rho = \overrightarrow{AA_\rho}$. Then
\begin{equation*}
  \Vb{x}_\rho = \Vb{b} + \rho \, \Vb{c}, \quad
 \Vb{y}_\rho = \Vb{a} + \rho\, \Vb{b} \quad
 \text{and} \quad
 \Vb{z}_\rho  = \Vb{c} + \rho \, \Vb{a}.
\end{equation*}
Since $\Vb{a}+ \Vb{b} + \Vb{c} = \Vb{0}$, we have
\begin{equation*}
\Vb{x}_\rho + \Vb{y}_\rho + \Vb{z}_\rho = \Vb{b} + \rho\, \Vb{c} + \Vb{a} + \rho\, \Vb{b}
+ \Vb{c} + \rho\, \Vb{a} =  (1+\rho) (\Vb{a} + \Vb{b} + \Vb{c}) = \Vb{0}.
\end{equation*}
Therefore, there exists an oriented triangle $XY\!Z$ whose sides have
the lengths  $x_\rho := Y\!Z = CC_\rho, y_\rho = ZX = BB_\rho$, and
$z_\rho = XY = AA_\rho$.  Here, as always in this paper, the vertices $X, Y, Z$ are labeled counterclockwise.

\begin{figure}[ht]

\psfrag{Z}[][]{\begin{picture}(0,0)
            \put(-2,-15){\makebox(0,0)[l]{$Z$}}
                        \end{picture}}

\psfrag{Y}[][]{\begin{picture}(0,0)
            \put(-3,0){\makebox(0,0)[l]{$Y$}}
                        \end{picture}}

\psfrag{X}[][]{\begin{picture}(0,0)
            \put(-8,-2){\makebox(0,0)[l]{$X$}}
                        \end{picture}}

\psfrag{B}[][]{\begin{picture}(0,0)
            \put(-1,-2){\makebox(0,0)[l]{$B$}}
                        \end{picture}}

\psfrag{Ap}[][]{\begin{picture}(0,0)
            \put(-6,-2){\makebox(0,0)[l]{$A_\rho$}}
                        \end{picture}}

\psfrag{Cp}[][]{\begin{picture}(0,0)
            \put(-20,4){\makebox(0,0)[l]{$C_\rho$}}
                             \end{picture}}

\psfrag{A}[][]{\begin{picture}(0,0)
            \put(-10,2){\makebox(0,0)[l]{$A$}}
                        \end{picture}}

\psfrag{Bp}[][]{\begin{picture}(0,0)
            \put(-5,2){\makebox(0,0)[l]{$B_\rho$}}
                        \end{picture}}

\psfrag{C}[][]{\begin{picture}(0,0)
            \put(-6,-2){\makebox(0,0)[l]{$C$}}
                        \end{picture}}

\psfrag{App}[][]{\begin{picture}(0,0)
            \put(-10,-12){\makebox(0,0)[l]{$A_\rho'$}}
                        \end{picture}}

\psfrag{D}[][]{\begin{picture}(0,0)
            \put(-8,2){\makebox(0,0)[l]{$D$}}
                        \end{picture}}

\setlength{\abovecaptionskip}{6pt}%
\setlength{\belowcaptionskip}{-6pt}%

\resizebox{!}{!}{%
  \includegraphics{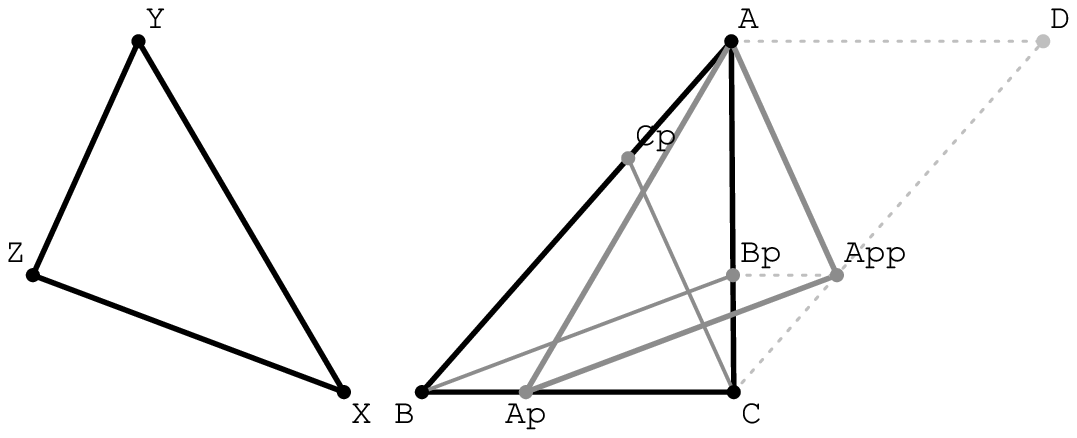}}
   \caption{Ceva's triangle ${\CT}_\rho (T)$ with $\rho =1/3$}
\label{f1}
\end{figure}

For $\rho\in (0, 1)$, there is a natural geometric construction of
the oriented triangle made by the three cevians which is worth recalling here since it is a straightforward modification of the one for the median
triangle. Let $D$ denote the point in the plane of $ABC$ such that
the quadrilateral $ABCD$ is a parallelogram having the diagonals
$AC$ and $BD$. The point $A_\rho$ on the segment $BC$ is such that
$BA_\rho = \rho BC$. Let $A_\rho'$ be the point on $C\!D$ such that
$C\!A_\rho' = \rho C\!D$. The sides of the triangle $A A_\rho A_\rho'$
are clearly equal to the three given cevians, see Figure~\ref{f1}.
We recognize the oriented triangle $XY\!Z$ as a reflection of the copy produced
by Hajja in \cite[Theorem~3.3]{Hajja1}.

Using the classic theorem known as Stewart's theorem (see for example \cite[Exercise~4 of Section~1.2]{Coxeter}), or simply applying the law of cosines, it is easy to calculate the lengths $x_\rho, y_\rho, z_\rho$:
 {\allowdisplaybreaks
\begin{align} \label{Crho}
x_\rho = CC_\rho
 & = \sqrt{\rho a^2 + (1-\rho)b^2  +  \rho(\rho - 1)c^2}, \\ \label{Brho} y_\rho = BB_\rho
  & = \sqrt{(1-\rho) a^2 + \rho(\rho-1)b^2 + \rho c^2}, \\ \label{Arho} z_\rho =  AA_\rho
   & = \sqrt{\rho(\rho-1)a^2  + \rho b^2 + (1-\rho)c^2}.
\end{align}}
\!\!Thus, starting with an ordered triple of sides $(a,b,c)$ and $\rho
\in \nR$, the ordered triple $(x_\rho,y_\rho,z_\rho)$ is uniquely
determined. We define Ceva's operator $\CT_\rho$ by
\[
\CT_\rho(a,b,c) : = (x_\rho, y_\rho, z_\rho).
\]
The oriented triangle $(x_\rho, y_\rho, z_\rho)$ we call {\em Ceva's  triangle of} $T$. The difference between our definition and the corresponding definition in \cite{Hajja1, Hajja2} is in the order of sides. There, the generalized median operator was defined by
\[
\mathcal H_\rho (a, b, c) = (z_\rho, y_\rho, x_\rho ).
\]
This innocent detail, however, creates problems in the iterative process investigated in \cite{Hajja1, Hajja2}. Our Ceva's operator behaves better  precisely due to the cyclical symmetry property observed in  \cite[Propositions~1~and~2]{Griffiths}. As we shall soon see, the operator $\CT_\rho$ produces a binary sequence of triangles, while ${\mathcal H}_\rho$ does not. Indeed, this is because ${\mathcal H}_\rho\!\circ {\mathcal H}_\rho = \CT_\rho \circ \CT_{1-\rho}$. Here, $\circ$ denotes composition of functions. We set $\CT_\rho^1 := \CT_\rho$ and recursively $\CT_\rho^{n+1} := \CT_\rho^n\circ \CT_\rho$ for all $n\in \nN$.

We often use capital letters $T, V, \ldots$ to denote oriented
triangles. Then $\CT_\rho(T)$, $\CT_\rho(V), \ldots$ denote
corresponding Ceva's triangles. We immediately note that, if $T$ is an equilateral triangle with side-length $a$, then $\CT_\rho (T)$ is also an equilateral triangle of side-length $a\sqrt{1-\rho+\rho^2}$. Because of this, the discussion below will only be concerned with non-equilateral
triangles.

Notice that \( \CT_0(a,b,c) = (b,a,c) \) and \( \CT_1(a,b,c) =
(a,c,b). \)  For completeness, we also define
 \(
\CT_\infty(a,b,c) = (c,b,a).
\)  As a consequence, the triangles $\CT_0(T), \CT_1(T)$ and
$\CT_\infty(T)$ are directly congruent to each other, and each is reversely congruent to $T$. We will see later that the
set ${\mathbb S} := \{0,1,\infty\}$ will play an important role
whenever we encounter direct similarity. Another important set is the unit interval $\nI:=[0, 1)$.

To summarize, we have defined Ceva's operator $\CT_\rho$ for
any $\rho \in \nR\cup \{\infty\}$. For a subset ${\mathbb J}$ of $\nRc$,  we
will write $\CT_{\mathbb J} (T)$ for the family $\{\CT_\rho (T): \rho\in {\mathbb J}\}$.

\section{The cone} \label{scone}

In this section, we show that a triple $(a,b,c)$ of positive numbers represents the side-lengths of a triangle if and only if $(a^2,b^2,c^2) \in \cQ$, where $\cQ$ is the interior in the first octant of the cone
\begin{equation*} 
x^2+y^2+z^2 - 2(xy+yz+zx) = 0.
\end{equation*}
That is,
\begin{equation*} 
\cQ = \left\{ \left[\!\!
\begin{array}{c}
x \\ y \\ z
\end{array}\! \!
\right] \, : \, x,y,z > 0, \   x^2+y^2+z^2 < 2(xy+yz+zx) \right\}.
\end{equation*}
This fact was already observed in \cite{Griffiths} in connection with  Heron's area formula. For completeness, we give a direct proof that
$|a-b| < c < a+b$ if and only if $\bigl[a^2 \ b^2 \ c^2\bigr]^\top \in \cQ$.
We have the following equivalences:
 {\allowdisplaybreaks
\begin{align} \nonumber
|a-b| < c < a+b \quad & \Leftrightarrow \quad a^2 + b^2 -2ab < c^2 < a^2+b^2 + 2ab  \\
 \nonumber
 \quad  & \Leftrightarrow \quad  |a^2 + b^2 - c^2 | < 2 a b  \\
   \nonumber
 \quad  & \Leftrightarrow \quad
   ( a^2 + b^2 - c^2 )^2  < 4 a^2 b^2  \\
  \label{eqcon}
 \quad  & \Leftrightarrow \quad
  a^4 + b^4 + c^4 < 2 a^2 b^2 + 2b^2c^2 + 2c^2a^2  \\
  \nonumber
 \quad  & \Leftrightarrow \quad  \bigl[a^2 \ b^2 \ c^2\bigr]^\top \in \cQ.
\end{align}}
The inequality in \eqref{eqcon} is further equivalent to
\begin{equation} \label{eqinQ}
2 ( a^4 + b^4 + c^4 )  <  ( a^2 + b^2 + c^2 )^2.
\end{equation}
Now, taking the square root of both sides, we adjust the last
inequality to look like an inequality for the dot product of two unit
vectors:
\begin{equation} \label{eqcga}
\frac{a^2\cdot 1 + b^2 \cdot 1 + c^2 \cdot 1}%
{\sqrt{a^4 + b^4 + c^4} \, \sqrt{3} }
 > \sqrt{\frac{2}{3}}.
\end{equation}
Denote by $\gamma_{_T}$ the angle between the vectors $[a^2 \ b^2 \
c^2 ]^\top$ and $[1 \ 1 \ 1]^\top$. Then the last inequality yields
that $\cos (\gamma_{_T}) > \sqrt{2/3}$. In other words, $(a,b,c)$
are the side-lengths of a triangle if and only if the vector  $[ a^2
\ b^2 \ c^2]^\top$ is inside the cone centered around the diagonal
$x=y=z$ and with the angle at the vertex equal to $\arccos
\sqrt{2/3} = \arctan (1/\sqrt{2})$. We will call the angle
$\gamma_{_T} \in \bigl[0,\arctan(1/\sqrt{2}) \bigr)$ the {\em
cone angle} of the triangle $T$.

\begin{figure}[ht]

\psfrag{P}[][]{\begin{picture}(0,0)
            \put(-10,-3){\makebox(0,0)[l]{$P$}}
                        \end{picture}}

\psfrag{oA}[][]{\begin{picture}(0,0)
            \put(-4,5){\makebox(0,0)[l]{$\omega_{_T}$}}
                        \end{picture}}

\psfrag{oB}[][]{\begin{picture}(0,0)
            \put(-11,-3){\makebox(0,0)[l]{$\omega_{_T}$}}
                        \end{picture}}

\psfrag{oC}[][]{\begin{picture}(0,0)
            \put(1,-5){\makebox(0,0)[l]{$\omega_{_T}$}}
                        \end{picture}}

\psfrag{B}[][]{\begin{picture}(0,0)
            \put(-3,-2){\makebox(0,0)[l]{$B$}}
                        \end{picture}}

\psfrag{A}[][]{\begin{picture}(0,0)
            \put(-12,2){\makebox(0,0)[l]{$A$}}
                        \end{picture}}

\psfrag{C}[][]{\begin{picture}(0,0)
            \put(-9,-2){\makebox(0,0)[l]{$C$}}
                        \end{picture}}

\setlength{\abovecaptionskip}{2pt}%
\setlength{\belowcaptionskip}{-6pt}%

\resizebox{!}{!}{%
  \includegraphics{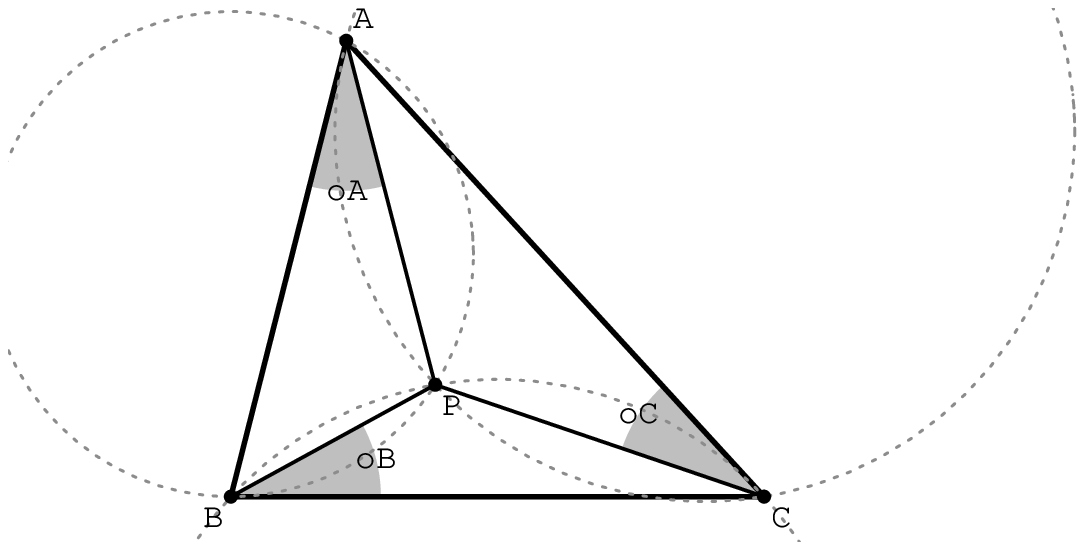}}
   \caption{The Brocard angle $\omega_{_T}$ of $T=ABC$}
\label{fb}
\end{figure}

In the next proposition we prove that the cone angle of $T$ uniquely determines another important angle of $T$, its Brocard angle, and vice versa. To define the Brocard angle of an oriented triangle $T = ABC$ one first proves that there exists a unique point $P$ such that the angles $PAB, PBC$ and $PCA$ (marked in gray in Figure~\ref{fb}) are equal to each other. This common angle is called the {\em Brocard angle} of $T$; it is denoted by $\omega_{_T}$. For more on this topic we refer to \cite[Chapters~XVI and~XVII]{Johnson} as a classical reference, or the more recent \cite[Chapter Ten]{Honsberger}.

\begin{proposition} \label{pbaca}
Let $T$ be a triangle, let  $\gamma_{_T}$ be its cone angle, and let $\omega_{_T}$ be its Brocard angle. Then
\begin{equation} \label{eqbaca}
3 (\tan \omega_{_T}\!)^2  + 2 (\tan \gamma_{_T}\!)^2 = 1.
 \end{equation}
Let $V$ also be a triangle. Then $\gamma_{_T} = \gamma_{_V}$ if and only if    $\omega_{_T} = \omega_{_V}$.
\end{proposition}
\begin{proof}
Let $T=(a,b,c)$. Following \cite[Theorem~2.4]{Hajja2}, we set
\[
k = \frac{a^4 + b^4 + c^4}{a^2b^2+b^2c^2+c^2a^2}.
\]
Using \eqref{eqcga} and the definition of $\gamma_{_T}$, we calculate $(\tan \gamma_{_T}\!)^2 = (2k-2)/(k+2)$. This
and the identity $(\tan \omega_{_T}\!)^2 = (2-k)/(k+2)$ from
\cite[Theorem~2.4]{Hajja2} yield \eqref{eqbaca}. Since by the  Cauchy-Schwarz inequality and \eqref{eqcon}, $1 \leq k < 2$, we have $0 < (\tan \omega_{_T}\!)^2 \leq 1/3$. Therefore, $\omega_{_T}, \omega_{_V} \in (0,\pi/6]$. As we already observed that $\gamma_{_T}, \gamma_{_V} \in \bigl[0,\arctan(1/\sqrt{2}) \bigr)$, the second claim in the proposition  follows from \eqref{eqbaca}.
\end{proof}

\section{Reflection matrices} \label{refm}

Let $\rho \in \nR$. Set $\pn{\rho} := \sqrt{1-\rho+\rho^2}$, $\pn{\infty} : = 1$,  and consider the left-circulant orthogonal matrices
\begin{equation} \label{eqMr}
\Vs{M}_\rho  = \frac{1}{\pn{\rho}^2} \left[\!\!
\begin{array}{ccc}
 \rho  & \ 1 - \rho & \ \rho(\rho-1) \\[8pt]
 1-\rho  &  \ \rho(\rho -1)  & \rho  \\[8pt]
 \rho(\rho-1)  & \rho  & 1-\rho
\end{array}\!\!
\right]
\  \text{and} \
\Vs{M}_\infty  = \left[\!\!
\begin{array}{ccc}
 0  & 0 & 1 \\[3pt]
 0  &  1  & 0  \\[3pt]
 1 & 0  & 0
\end{array}\!\!
\right].
\end{equation}
We note that it can be shown that $\bigl\{\Vs{M}_\rho, -\Vs{M}_\rho:\, \rho\in \nR\cup\{\infty\}\bigr\}$ is the family of all left-circulant orthogonal $3\times 3$ matrices.

It follows from \eqref{Crho},\eqref{Brho},\eqref{Arho} that the squares of the side-lengths of Ceva's triangle $\CT_\rho(T)$ are related to the squares of the side-lengths of the original oriented triangle $T$ in the following simple way:
\begin{equation}
\label{eqsqMsq}
\left[\!\!\begin{array}{c}
 x_\rho^2 \\[5pt] y_\rho^2 \\[5pt] z_\rho^2 \end{array}\!\!\right]
  =
 \pn{\rho}^2 \Vs{M}_\rho \left[\!\!\begin{array}{c}
 a^2 \\[5pt] b^2 \\[5pt] c^2 \end{array}\!\!\right], \qquad \rho \in \nRc.
\end{equation}

The fact that for every triangle $T$, $\CT_\rho(T)$ is also a triangle is  equivalent to the statement that the matrix $\Vs{M}_\rho$ maps $\cQ$ into $\cQ$. We proved this geometrically at the beginning of the paper and it is proved as a matrix statement in \cite[Propositions~1 and~2]{Griffiths}. However, to fully understand the family of triangles $\CT_\rho(T), \rho \in \nRc$, we need a deeper understanding of the family of matrices in \eqref{eqMr}. This and the following three sections provide that understanding.

For an arbitrary $\rho \in \nRc$, the matrix $\Vs{M}_\rho$ is symmetric and orthogonal. Hence its eigenvalues are $1$ and $-1$ and there is an orthonormal basis consisting of eigenvectors of $\Vs{M}_\rho$. To find such a basis, we first observe that the row sums of each $\Vs{M}_\rho$ are equal to $1$ making the vector $[1\ 1 \ 1]^\top$ an eigenvector corresponding to the eigenvalue $1$. We normalize the opposite of this vector and calculate the orthonormal positively oriented eigenvectors of $\Vs{M}_\rho$ corresponding to the eigenvalues $1$, $-1$, $1$, respectively, to be
\begin{alignat*}{3} 
\Vb{p}_\rho &:=  \frac{1}{\sqrt{6}\ \pn{\rho}}\left[\!\!\begin{array}{c}
 1+\rho \\[5pt] 1- 2\rho \\[3pt] \rho -2 \end{array}\!\!\right],
  &  \qquad
  \Vb{q}_\rho &:= \frac{1}{\sqrt{2}\ \pn{\rho}}\left[\!\!\begin{array}{c}
 1-\rho \\[5pt] -1 \\[3pt] \rho \end{array}\!\!\right],
  &  \qquad
  \Vb{r} & := \frac{-1}{\sqrt{3}} \left[\!\begin{array}{r}
  1 \\ 1 \\ 1 \end{array}\!\right]. \\
\intertext{The corresponding eigenvectors of $\Vs{M}_\infty$ are}
\Vb{p}_\infty & :=  \frac{1}{\sqrt{6}} \left[\!\!\begin{array}{r}
  -1 \\ 2 \\ -1 \end{array}\!\right],  &
 \Vb{q}_\infty & :=\frac{1}{\sqrt{2}} \left[\!\!\begin{array}{r}
 1 \\ 0 \\ -1 \end{array}\!\right], &
 \Vb{r} & := \frac{-1}{\sqrt{3}} \left[\!\begin{array}{r}
  1 \\ 1 \\ 1 \end{array}\!\right].
\end{alignat*}
Consequently, the matrix $\Vs{M}_\rho, \rho \in \nRc$, induces the reflection with respect to the plane spanned by the vectors $\Vb{p}_\rho$ and $\Vb{r}$. Thus, $\Vs{M}_\rho$ is a reflection matrix.

\begin{remark} \label{risop}
The reflection planes corresponding to $\Vs{M}_\infty, \Vs{M}_0$ and $\Vs{M}_1$ are given by the equations $x=z$, $x=y$, and $y=z$, respectively. Therefore, the triples in the intersection of these planes with $\cQ$ correspond to isosceles triangles. Moreover, the triples that are in $\cQ$ and in the quadrants determined by $-\Vb{r}$ and each of $\Vb{p}_\infty, \Vb{p}_0, \Vb{p}_1$ correspond to wide triangles and the triples that are in $\cQ$  and in the quadrants determined by $-\Vb{r}$ and each of $-\Vb{p}_\infty, -\Vb{p}_0, -\Vb{p}_1$ are narrow.
\end{remark}

Next, we will prove that an arbitrary reflection across a plane which contains the vector $\Vb{r}$ is in the family \eqref{eqMr}. Such a plane is uniquely determined by its trace in the plane spanned by the vectors
\(\Vb{p}_0, \Vb{q}_0.\) In turn, this trace is uniquely determined by its angle $\vartheta \in \bigl[-\tfrac{\pi}{3},\tfrac{2\pi}{3}\bigr)$ with the vector $\Vb{p}_0$.

\begin{figure}

\setlength{\abovecaptionskip}{2pt}%
\setlength{\belowcaptionskip}{-8pt}%

\!\!\begin{minipage}{0.64\textwidth}

\psfrag{op3}[][]{\begin{picture}(0,0)
            \put(16,0){\makebox(0,0)[l]{$-\tfrac{\pi}{3}$}}
                        \end{picture}}
\psfrag{op6}[][]{\begin{picture}(0,0)
            \put(16,0){\makebox(0,0)[l]{$-\tfrac{\pi}{6}$}}
                        \end{picture}}
\psfrag{p3}[][]{\begin{picture}(0,0)
            \put(0,0){\makebox(0,0)[l]{$\tfrac{\pi}{3}$}}
                        \end{picture}}
\psfrag{p6}[][]{\begin{picture}(0,0)
            \put(0,0){\makebox(0,0)[l]{$\tfrac{\pi}{6}$}}
                        \end{picture}}
\psfrag{p2}[][]{\begin{picture}(0,0)
            \put(0,0){\makebox(0,0)[l]{$\tfrac{\pi}{2}$}}
                        \end{picture}}
\psfrag{p23}[][]{\begin{picture}(0,0)
            \put(0,0){\makebox(0,0)[l]{$\tfrac{2\pi}{3}$}}
                        \end{picture}}

\psfrag{ot}[][]{\begin{picture}(0,0)
            \put(-7,18){\makebox(0,0)[l]{$-2$}}
                        \end{picture}}
\psfrag{oo}[][]{\begin{picture}(0,0)
            \put(-7,18){\makebox(0,0)[l]{$-1$}}
                        \end{picture}}
\psfrag{h}[][]{\begin{picture}(0,0)
            \put(-2,-1){\makebox(0,0)[l]{$\tfrac{1}{2}$}}
                        \end{picture}}
\psfrag{o}[][]{\begin{picture}(0,0)
            \put(-3,2){\makebox(0,0)[l]{$1$}}
                        \end{picture}}
\psfrag{t}[][]{\begin{picture}(0,0)
            \put(-3,2){\makebox(0,0)[l]{$2$}}
                        \end{picture}}
\resizebox{!}{!}{%
  \includegraphics{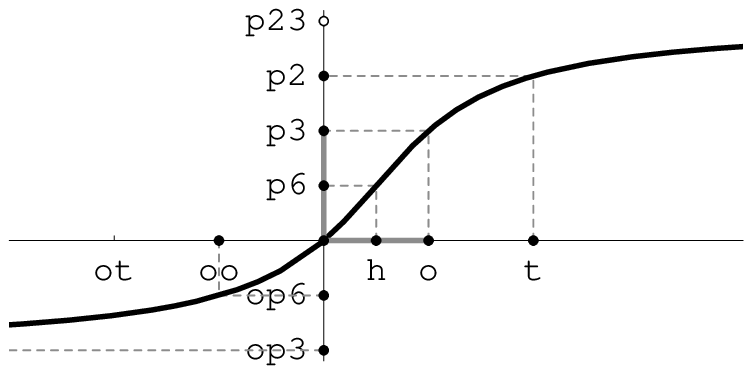}}
    \caption{The function $\Phi^{-1}$}
\label{fPhii}
\end{minipage}   
\begin{minipage}{0.34\textwidth}

\vspace*{10pt}

\psfrag{p3}[][]{\begin{picture}(0,0)
            \put(-2,0){\makebox(0,0)[l]{$\tfrac{\pi}{3}$}}
                        \end{picture}}
\psfrag{p6}[][]{\begin{picture}(0,0)
            \put(-2,0){\makebox(0,0)[l]{$\tfrac{\pi}{6}$}}
                        \end{picture}}

\psfrag{zx}[][]{\begin{picture}(0,0)
            \put(-10,0){\makebox(0,0)[l]{$0$}}
                        \end{picture}}

\psfrag{h}[][]{\begin{picture}(0,0)
            \put(-3,-3){\makebox(0,0)[l]{$\tfrac{1}{2}$}}
                        \end{picture}}
\psfrag{o}[][]{\begin{picture}(0,0)
            \put(-3,-1){\makebox(0,0)[l]{$1$}}
                        \end{picture}}

\setlength{\abovecaptionskip}{5pt}%
\setlength{\belowcaptionskip}{-4pt}%

\rule{5pt}{0pt} \resizebox{!}{!}{%
  \includegraphics{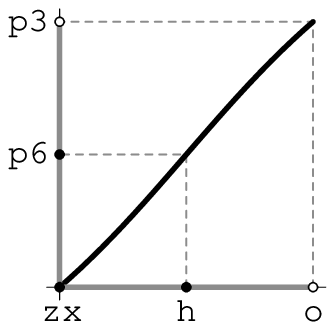}}
    \caption{$\Phi^{-1}$ on $\nI$}
\label{fPhiiI}
\end{minipage}

\end{figure}

Denote by $\vartheta_\rho$ the angle between $\Vb{p}_\rho$ and $\Vb{p}_0$. Then,
\[
\cos \vartheta_\rho = \Vb{p}_\rho \cdot \Vb{p}_0
  = \frac{2-\rho}{2\, \pn{\rho}}, \quad
 \sin \vartheta_\rho = \Vb{p}_\rho \cdot \Vb{q}_0
   = \frac{\sqrt{3} \rho}{2\, \pn{\rho}},
 \quad \text{and} \quad
\tan \vartheta_\rho = \frac{\sqrt{3} \rho}{2-\rho}.
\]
Solving the last equation for $\rho$ we get
\[
\rho = \dfrac{2 \, \tan \vartheta_\rho}{\sqrt{3}+\tan \vartheta_\rho} = \frac{\sqrt{3}}{2} \tan\Bigl(\vartheta_\rho-\frac{\pi}{6}\Bigr)+\frac{1}{2} = \frac{\sin \vartheta_\rho}{\cos\bigl(\vartheta_\rho - \tfrac{\pi}{6}\bigr)}.
\]
We define now
\begin{equation*} 
\Phi(\vartheta) := \left\{
\begin{array}{cl}
 \infty
  & \ \text{if} \quad \vartheta = -\tfrac{\pi}{3}, \\[6pt]
\dfrac{\sin \vartheta}{\cos \bigl(\vartheta - \tfrac{\pi}{6}\bigr)}
   & \ \text{if} \quad \vartheta \in \bigl(-\tfrac{\pi}{3},\tfrac{2\pi}{3}\bigr).
\end{array}
 \right.
\end{equation*}
Clearly, $\Phi$ is an increasing bijection between  $\bigl(-\tfrac{\pi}{3},\tfrac{2\pi}{3}\bigr)$ and $\nR$. The inverse of the function $\Phi$ is, see Figures~\ref{fPhii} and~\ref{fPhiiI},
\[
\Phi^{-1}(\rho) = \left\{
\begin{array}{cl}
  - \tfrac{\pi}{3}
  & \ \text{if} \quad \rho = \infty,  \\[10pt]
 \arctan\Bigl(\tfrac{2}{\sqrt{3}}\bigl(\rho - \tfrac{1}{2} \bigr)\Bigr) +\tfrac{\pi}{6}
   & \ \text{if} \quad \rho\in\nR.
\end{array}
 \right.
\]

Since the function $\Phi$ is a bijection between $\bigl[-\tfrac{\pi}{3},\tfrac{2\pi}{3}\bigr)$ and $\nRc$, that is since the range of $\Phi^{-1}$ is the interval $\bigl[-\tfrac{\pi}{3},\tfrac{2\pi}{3}\bigr)$, all the reflections across planes containing the vector $\Vb{r}$ are represented in the family \eqref{eqMr}.

\section{Three special groups}

The interval $\bigl[-\tfrac{\pi}{3},\tfrac{2\pi}{3}\bigr)$ with the addition modulo $\pi$, which we denote by $\oplus$, is a commutative group. The identity element is $0$. The inverse of $\vartheta \in \bigl[-\tfrac{\pi}{3},\tfrac{\pi}{3}\bigr]$ is $-\vartheta$ and the inverse of $\vartheta \in \bigl(\tfrac{\pi}{3},\tfrac{2\pi}{3}\bigr)$ is $\pi - \vartheta$. The bijection $\Phi$ then induces a natural group structure on $\nRc$. We denote by $\Bo$ the operation of this group. We have
\[
\rho\Bo \tau := \Phi\bigl(\Phi^{-1}(\rho)+ \Phi^{-1}(\tau)\bigr), \qquad \rho,\tau \in \nRc.
\]
Let $\rho, \tau \in \nR\setminus\{2\}$ be such that $\rho\tau \neq 1$. Set $\vartheta_\rho = \Phi^{-1}(\rho), \vartheta_\tau = \Phi^{-1}(\tau)$. Then,
{\allowdisplaybreaks
\begin{align*}
\rho\Bo \tau & = \Phi(\vartheta_\rho + \vartheta_\tau) \\
 & = \dfrac{2\tan(\vartheta_\rho+\vartheta_\tau)}%
      {\sqrt{3}+\tan(\vartheta_\rho+\vartheta_\tau)} \\
& = \dfrac{2\dfrac{\tan\vartheta_\rho+\tan\vartheta_\tau}%
           {1-(\tan\vartheta_\rho)(\tan \vartheta_\tau)}}%
      {\sqrt{3}+\dfrac{\tan \vartheta_\rho +\tan \vartheta_\tau}%
           {1-(\tan\vartheta_\rho)(\tan\vartheta_\tau)}} \\
& = \dfrac{2 \bigl(\tan \vartheta_\rho +\tan \vartheta_\tau \bigr)}%
      {\sqrt{3}\bigl(1-(\tan\vartheta_\rho)(\tan\vartheta_\tau)\bigr)+ \tan \vartheta_\rho +\tan \vartheta_\tau} \\
 & = \dfrac{2 \frac{\sqrt{3} \rho}{2-\rho}
     +2 \frac{\sqrt{3} \tau}{2-\tau}}%
      {\sqrt{3}
       -\sqrt{3}\frac{\sqrt{3}\rho}{2-\rho}\frac{\sqrt{3}\tau}{2-\tau}+ \frac{\sqrt{3}\rho}{2-\rho}+\frac{\sqrt{3}\tau}{2-\tau}} \\
 & = \frac{2\rho(2-\tau) + 2\tau(2-\rho)}%
    {(2-\rho)(2-\tau)-3\rho\tau+\rho(2-\tau) +\tau(2-\rho)}\\
     & = \frac{4\rho + 4\tau -4 \rho\tau}%
    {4-2\rho-2\tau+\rho\tau-3\rho\tau+2\rho + 2\tau -2 \rho\tau}\\
    & = \frac{\rho + \tau -\rho\tau}%
    {1-\rho\tau}.
\end{align*}
}
The other values of $\rho, \tau \in \nRc$ are treated similarly to get
\[
\rho \Bo \tau  =  \left\{ \begin{array}{cl}
 \frac{\rho + \tau - \rho\tau}{1- \rho \tau} \quad & \text{if} \quad \rho,\tau \in \nR, \ \ \rho\tau \neq 1, \\[4pt]
 \infty \quad & \text{if} \quad \rho,\tau \in \nR, \ \ \rho\tau = 1, \\[4pt]
 1 - \frac{1}{\rho} \quad & \text{if} \quad \rho \in \nR\setminus\{0\}, \ \ \tau = \infty, \\[4pt]
 \infty \quad & \text{if} \quad \rho = 0, \ \tau = \infty \ \text{or} \ \rho =\infty, \ \tau = 0, \\[4pt]
  1 - \frac{1}{\tau} \quad & \text{if} \quad  \rho = \infty, \ \ \tau \in \nR\setminus\{0\}, \\[4pt]
  1 \quad & \text{if} \quad \rho = \infty, \ \ \tau = \infty.
\end{array} \right.
\]

The set $\nRc$ with the operation $\Bo$ is a commutative group with the identity element $0$. The inverses are
\[
\Ton{\rho} := \left\{\begin{array}{cl}
\frac{\rho}{\rho-1}  & \ \text{if} \quad \rho \in \nR\setminus\{1\}, \\
\infty   & \ \text{if} \quad \rho = 1, \\
1      & \ \text{if} \quad \rho = \infty.
\end{array}\right.
\]
The set $\nS := \{0,1,\infty\}$ is a cyclic subgroup of $(\nRc,\Bo)$
of order $3$ which corresponds to the cyclic subgroup
$\bigl\{0,\frac{\pi}{3},-\frac{\pi}{3}\bigr\}$ of
$\Bigl(\bigl[-\frac{\pi}{3},\frac{2\pi}{3}\bigr),\oplus\Bigl)$.
Similarly, $\nT := \{0, 1/2, 1, 2, \infty, -1\}$ is a cyclic
subgroup of $(\nRc,\Bo)$ of order $6$ which corresponds to the
cyclic subgroup
$\bigl\{0,\frac{\pi}{6},\frac{\pi}{3},\frac{\pi}{2},-\frac{\pi}{3},-\frac{\pi}{6}\bigr\}$
of $\Bigl(\bigl[-\frac{\pi}{3},\frac{2\pi}{3}\bigr),\oplus\Bigl)$.
The $\Bo$-operation on $\nT$ is summarized in Table~\ref{ta1}.


For the three special values of $\tau \in \nS$, the operation $\Bo$ gives three functions that we will encounter in the definition of the function $\pI$ below:
\begin{equation*}
 0 \Bo \rho = \rho, \quad  1 \Bo \rho  = \frac{1}{1-\rho}, \quad \infty \Bo \rho  = 1 - \frac{1}{\rho}, \quad \rho \in \nRc.
\end{equation*}
We will write $\nS\Bo \rho$ for the set $\{0 \Bo \rho, 1 \Bo \rho, \infty\Bo \rho\}$. Another interesting set of functions is $\nS\Bo \Ton{\rho} = \nS\Bo \tfrac{1}{\rho}$:
\[
0 \Bo \Ton{\rho} = 1 \Bo \frac{1}{\rho} = \frac{\rho}{\rho-1},
 \quad
  1 \Bo \Ton{\rho}  = \infty \Bo \frac{1}{\rho} = 1-\rho,
  \quad
  \infty \Bo \Ton{\rho} = 0 \Bo \frac{1}{\rho} = \frac{1}{\rho},
  \quad
 \rho \in \nRc.
\]

\begin{table}[ht]

\setlength{\abovecaptionskip}{-6pt}%
\setlength{\belowcaptionskip}{8pt}%

\begin{center}
     \renewcommand*\arraystretch{1.5}
\begin{tabular}{@{\extracolsep{0pt}}|c"c|c|c|c|c|c|}\hline
 $\Bo$ & $0$ & $1/2$ & $1$ &  $2$ & $\infty$ & $-1$ \\\thickhline
$0$      &$0$     &$1/2$   &$1$     &$2$     &$\infty$ & $-1$ \\\hline
$1/2$    &$1/2$   &$1$     &$2$     &$\infty$&$-1$    & $0$ \\\hline
$1$      &$1$     &$2$     &$\infty$&$-1$    &$0$     & $1/2$ \\\hline
$2$      &$2$     &$\infty$&$-1$    &$0$     &$1/2$    & $1$ \\\hline
$\infty$ &$\infty$&$-1$    &$0$     &$1/2$   &$1$      &$2$ \\\hline
$-1$     &$-1$   &$0$      &$1/2$   &$1$     &$2$     &$\infty$ \\\hline
\end{tabular}
\caption{$\Bo$-operation table on $\nT$} \label{ta1}
\end{center}

\end{table}

We remark that the set $\nS\Bo\bigl\{\rho, \tfrac{1}{\rho} \bigr\}=\bigl(\nS\Bo\rho\bigr)\cup \bigl(\nS\Bo\tfrac{1}{\rho}\bigr)$, which consists of the six special functions above, represents all  possible values of the cross-ratio of four numbers. Furthermore, in this
context, the set $\nT$ consists of the only extended real-valued
fixed points of a cross-ratio. Therefore, it is plausible that there
exists a further connection between the group structure $\Bo$ and
cross-ratios. However, since it is unclear to us whether this link
would present any further simplifications in our work, we do not
pursue it further here.

\begin{figure}

\psfrag{yz}[][]{\begin{picture}(0,0)
            \put(3,2){\makebox(0,0)[l]{$0$}}
                        \end{picture}}
\psfrag{yo}[][]{\begin{picture}(0,0)
            \put(3,2){\makebox(0,0)[l]{$1$}}
                        \end{picture}}
\psfrag{ot}[][]{\begin{picture}(0,0)
            \put(-7,2){\makebox(0,0)[l]{$-2$}}
                        \end{picture}}
\psfrag{oo}[][]{\begin{picture}(0,0)
            \put(-7,2){\makebox(0,0)[l]{$-1$}}
                        \end{picture}}
\psfrag{z}[][]{\begin{picture}(0,0)
            \put(-3,2){\makebox(0,0)[l]{$0$}}
                        \end{picture}}
\psfrag{h}[][]{\begin{picture}(0,0)
            \put(-3,1){\makebox(0,0)[l]{$\tfrac{1}{2}$}}
                        \end{picture}}
\psfrag{o}[][]{\begin{picture}(0,0)
            \put(-3,2){\makebox(0,0)[l]{$1$}}
                        \end{picture}}
\psfrag{t}[][]{\begin{picture}(0,0)
            \put(-3,2){\makebox(0,0)[l]{$2$}}
                        \end{picture}}
\psfrag{th}[][]{\begin{picture}(0,0)
            \put(-3,2){\makebox(0,0)[l]{$3$}}
                        \end{picture}}

\setlength{\abovecaptionskip}{2pt}%
\setlength{\belowcaptionskip}{-8pt}%

\hspace*{-10pt}\resizebox{!}{!}{%
  \includegraphics{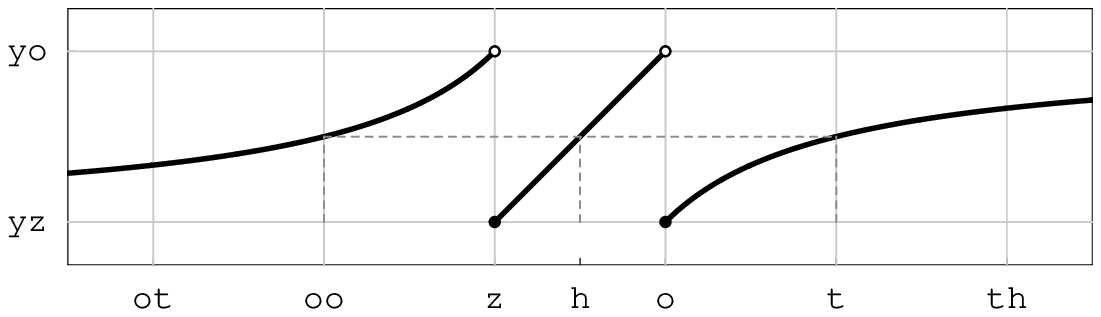}}
    \caption{The function $\pI$}
\label{funf}
\end{figure}

It turns out that the factor group $(\nRc)/\nS$ plays an important role in this paper. This factor group is isomorphic to the factor group $\bigl[\!-\tfrac{\pi}{3},\tfrac{2\pi}{3}\!\bigr)/\bigl\{0,\tfrac{\pi}{3},-\tfrac{\pi}{3}\bigr\}$, which in turn is isomorphic to the group $\bigl[0,\tfrac{\pi}{3}\bigr)$ with the addition modulo $\tfrac{\pi}{3}$. Since $\Phi\Bigl(\bigl[0,\tfrac{\pi}{3}\bigr) \Bigr) = \nI=[0, 1)$, $\Phi$ induces a natural group structure on $\nI$. A different way of understanding this group on $\nI$ is to notice that, for an arbitrary $\rho \in \nRc$, the intersections of the sets $\nS\Bo \rho$ and $\nI$ consist of exactly one number, which we denote by $\pI(\rho)$. That is, see Figure~\ref{funf},
\begin{equation*} 
\pI(\rho) = \left\{\!\! \begin{array}{rl}
1 \Bo \rho  & \text{if} \quad \rho = \infty \ \text{or} \ \rho < 0, \\
0 \Bo \rho  & \text{if} \quad \rho \in \nI, \\
\infty \Bo \rho & \text{if} \quad \rho \geq 1.
\end{array}\right.
\end{equation*}
The induced group operation on $\nI$ is denoted by $\BoI$:
\[
\rho \BoI \tau := \pI\bigl(\rho \Bo \tau\bigr), \qquad \rho, \tau \in \nI.
\]
With this definition, $\pI: \nRc \to \nI$ is an endomorphism between the groups $(\nRc,\Bo)$ and $(\nI,\BoI)$.  For $\rho, \tau \in \nI$, we have $\rho\Bo \tau \geq 0$ and also $\rho + \tau < 1$ if and only if $\rho\Bo \tau < 1$. Therefore, the above definition is equivalent to
\[
\rho \BoI \tau = \left\{\begin{array}{rll}
\rho\Bo \tau  = &\!\!\!\!\! \frac{\rho+\tau-\rho\tau}{1-\rho\tau} & \ \text{if} \quad \rho, \tau \in \nI \ \ \text{and} \ \ \rho+\tau < 1, \\[6pt]
\infty \Bo \rho\Bo \tau  = &\!\!\!\!\!  \frac{\rho+\tau-1}{\rho+\tau-\rho\tau} & \ \text{if} \quad \rho, \tau \in \nI \ \ \text{and} \ \ \rho+\tau \geq 1.
\end{array} \right.
\]
Given $\rho \in \nI$, its inverse $\TonI{\rho}$ with respect to $\BoI$ is $1-\rho$, or, equivalently, $\TonI{\rho} = 1\Bo\Ton{\rho}$.

\begin{remark} \label{rappI}
It is interesting to compare the groups $(\nI,\BoI)$ and $(\nI,\oplus)$, where $\oplus$ denotes the addition modulo $1$. As noticed above $\TonI{\rho} = \ominus \rho$ for all $\rho \in \nI$. Here, as usual in additive groups, we denote opposite elements by using the notation $\ominus$. However, these two groups are different. It turns out that the maximum value of $|\rho \BoI \tau - \rho \oplus \tau|$ is less than $0.042$, while the maximum value of the relative error $\bigl| \bigl(\rho \BoI \tau - \rho \oplus \tau\bigr)/(\rho \BoI \tau)\bigr|$ is $7-4\sqrt{3} \approx 0.072$. However, since the groups $(\nI,\oplus)$ and $\Bigl(\bigl[0,\tfrac{\pi}{3}\bigr),\oplus\Bigr)$  are isomorphic (where $\oplus$ denotes the addition modulo the length of interval), the groups $(\nI,\BoI)$ and $(\nI,\oplus)$ are isomorphic with the isomorphism $\tfrac{3}{\pi}\Phi^{-1} : \nI \to \nI$.
\end{remark}

\section{Functions on groups}

Let $\fvk \in \nRc$ and consider the equation $\xi\Bo \xi = \fvk$. One can easily verify that this equation has a unique solution in $[-1,1)$, which we denote by $\sqrt[\Bo]{\fvk}$\,; it is given by
\begin{equation*} 
\sqrt[\Bo]{\fvk} := \left\{ \begin{array}{cc}
 \dfrac{\fvk}{1+\llparenthesis \fvk\rrparenthesis}
  & \quad \text{if} \quad \fvk\in \nR, \\[14pt]
 - 1 \ & \quad \text{if} \quad \fvk = \infty.
\end{array} \right.
\end{equation*}
In particular, $\sqrt[\Bo]{0} = 0, \sqrt[\Bo]{1} = 1/2$. The solution set of $\xi\Bo \xi = \fvk$ is $\sqrt[\Bo]{\fvk}\Bo \{0,2\}$.

If $\fvk \in \nI$, then $0 \leq \sqrt[\Bo]{\fvk} < \tfrac{1}{2}$ and $\sqrt[\Bo]{\fvk} < \fvk$. Therefore $\sqrt[\Bo]{\fvk}$ is a solution of the equation $\xi\BoI \xi = \fvk, \, \xi\in\nI$. To find the second solution of this equation recall that $\TonI{(\tfrac{1}{2})} = \tfrac{1}{2}$, and therefore
\[
\fvk\BoI\TonI{\bigl(\tfrac{1}{2}\BoI\sqrt[\Bo]{\fvk}\bigr)} =
\fvk\BoI\TonI{\bigl(\sqrt[\Bo]{\fvk}\bigr)} \BoI \tfrac{1}{2} = \tfrac{1}{2}\BoI \sqrt[\Bo]{\fvk}.
\]
Hence, the other solution is $\tfrac{1}{2}\BoI\sqrt[\Bo]{\fvk}$. Since $\tfrac{1}{2}+ \sqrt[\Bo]{\fvk} < 1$, we have $\tfrac{1}{2}\BoI \sqrt[\Bo]{\fvk} = \tfrac{1}{2}\Bo \sqrt[\Bo]{\fvk}$. As $\tfrac{1}{2}\Bo \tau \geq \tfrac{1}{2}$ for all $\tau \in \bigl[ 0,\tfrac{1}{2}\bigr)$, we have
\begin{equation} \label{eqfps}
0 \leq \sqrt[\Bo]{\fvk} < \tfrac{1}{2} \leq \tfrac{1}{2}\Bo \sqrt[\Bo]{\fvk} < 1.
\end{equation}
It is trivial to see that $\pI(\!\sqrt[\Bo]{\fvk}) = \sqrt[\Bo]{\pI(\fvk)}$ for $\fvk \in \nI$, but $\pI(\!\sqrt[\Bo]{\fvk}) =\tfrac{1}{2}\Bo\sqrt[\Bo]{\pI(\fvk)}$ for $\fvk \in (\nRc)\setminus \nI$.

Now consider the group $\Bigl(\bigl[0,\tfrac{\pi}{3}\bigr),\oplus \Bigr)$. As before, $\oplus$ stands for
the addition modulo $\tfrac{\pi}{3}$. Let $\varphi \in \bigl[0,\tfrac{\pi}{3}\bigr)$.  The solutions
of the equation $\vartheta \oplus \vartheta = \varphi, \, 0 \leq \vartheta < \tfrac{\pi}{3}$, are $\tfrac{\varphi}{2}$ and
$\tfrac{\varphi}{2}\oplus\tfrac{\pi}{6}$. Since $\Phi: \bigl[0,\tfrac{\pi}{3}\bigr) \to [0,1)$ is an increasing
isomorphism between the groups $\Bigl(\bigl[0,\tfrac{\pi}{3}\bigr),\oplus \Bigr)$ and $(\nI ,\BoI)$, we have
\[
\Phi\bigl(\tfrac{\varphi}{2}\bigr) = \sqrt[\Bo]{\Phi(\varphi)} \qquad \text{and} \qquad \Phi\bigl(\tfrac{\varphi}{2}
\oplus\tfrac{\pi}{6}\bigr) = \tfrac{1}{2}\Bo \sqrt[\Bo]{\Phi(\varphi)}.
\]

 \begin{figure}[ht]


\setlength{\abovecaptionskip}{6pt}%
\setlength{\belowcaptionskip}{-8pt}%

\begin{minipage}{0.49\textwidth}

\psfrag{rph}[][]{\begin{picture}(0,0)
     \put(-12,0){\makebox(0,0)[l]{$\tfrac{\pi}{3}-\varphi$}}
                        \end{picture}}

\psfrag{zx}[][]{\begin{picture}(0,0)
            \put(2,2){\makebox(0,0)[l]{$0$}}
                        \end{picture}}

\psfrag{px3}[][]{\begin{picture}(0,0)
            \put(-4,0){\makebox(0,0)[l]{$\tfrac{\pi}{3}$}}
                        \end{picture}}
\psfrag{px6}[][]{\begin{picture}(0,0)
            \put(-4,0){\makebox(0,0)[l]{$\tfrac{\pi}{6}$}}
                        \end{picture}}
\psfrag{py3}[][]{\begin{picture}(0,0)
            \put(0,-1){\makebox(0,0)[l]{$\tfrac{\pi}{3}$}}
                        \end{picture}}
\psfrag{py6}[][]{\begin{picture}(0,0)
            \put(0,0){\makebox(0,0)[l]{$\tfrac{\pi}{6}$}}
                        \end{picture}}

\resizebox{!}{!}{%
  \includegraphics{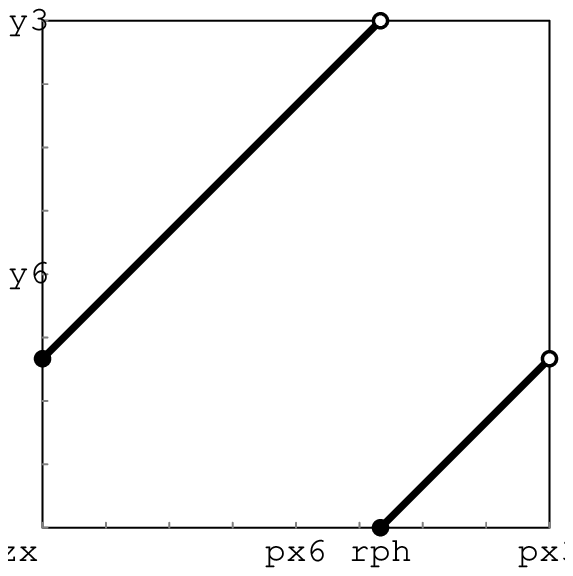}}
    \caption{$f_\varphi$ with $\varphi = \pi/9$}
\label{funFPi}
\end{minipage} \hfill
\begin{minipage}{0.49\textwidth}

\psfrag{ral}[][]{\begin{picture}(0,0)
            \put(-12,0){\makebox(0,0)[l]{$1-\fvk$}}
                        \end{picture}}
\psfrag{zx1}[][]{\begin{picture}(0,0)
            \put(2,2){\makebox(0,0)[l]{$0$}}
                        \end{picture}}

\psfrag{x1}[][]{\begin{picture}(0,0)
            \put(-4,0){\makebox(0,0)[l]{$1$}}
                        \end{picture}}
\psfrag{xh}[][]{\begin{picture}(0,0)
            \put(-4,-1){\makebox(0,0)[l]{$\tfrac{1}{2}$}}
                        \end{picture}}
\psfrag{y1}[][]{\begin{picture}(0,0)
            \put(0,0){\makebox(0,0)[l]{$1$}}
                        \end{picture}}
\psfrag{yh}[][]{\begin{picture}(0,0)
            \put(0,0){\makebox(0,0)[l]{$\tfrac{1}{2}$}}
                        \end{picture}}

\resizebox{!}{!}{%
  \includegraphics{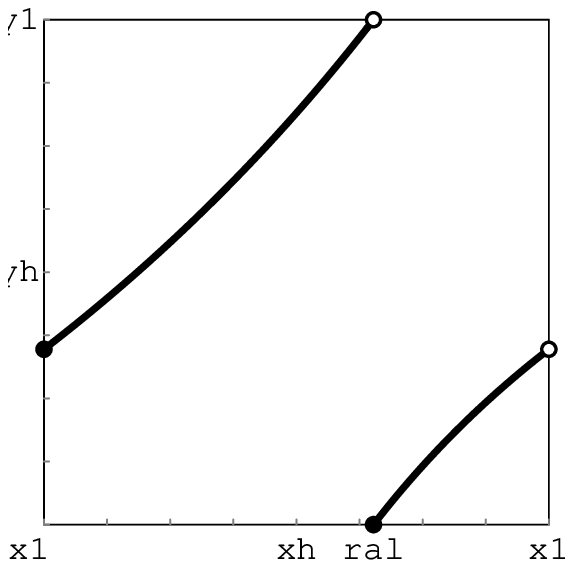}}
    \caption{\!\!\!$F_\fvk$,\! $\fvk=\Phi^{-1}(\tfrac{\pi}{9})\approx 0.35$}
\label{funFI}
\end{minipage}

\end{figure}

\begin{proposition} \label{padx}
Let \(\fvk, \lep, \rep  \in \nI\) and \(\lep < \rep\). Consider the function
\[
F_\fvk(\xi) : = \fvk\BoI \xi, \quad \xi\in\nI.
\]
\begin{enumerate}[{\rm (a)}]
\item  \label{padxia}
The function $F_\fvk$ is a bijection on $\nI$. Its inverse is $F_{\TonI{\fvk}}$.
\item \label{padxib}
If $F_\fvk(\lep) < F_\fvk(\rep)$, then $F_\fvk$ maps \([\lep, \rep]\) onto \([F_\fvk(\lep), F_\fvk(\rep)]\) as an increasing bijection.
\end{enumerate}
\end{proposition}
\begin{proof}
The statement (\ref{padxia}) is trivial. To prove (\ref{padxib}), let $\varphi \in \bigl(0,\tfrac{\pi}{3}\bigr)$ and consider the function $f_\varphi(\vartheta) = \varphi\oplus \vartheta$, see Figure~\ref{funFPi}. Since
\[
f_\varphi(\vartheta) = \left\{ \begin{array}{llrl}
 \varphi + \vartheta & \
    \text{if} & \quad 0  & \!\!\!\! \leq   \vartheta < \tfrac{\pi}{3} - \varphi, \\[4pt]
 \varphi + \vartheta - \tfrac{\pi}{3} & \
 \text{if} & \quad \tfrac{\pi}{3} - \varphi & \!\!\!\! \leq  \vartheta < \tfrac{\pi}{3},
 \end{array}
\right.
\]
we have the equivalence
\[
\vartheta_1 < \tfrac{\pi}{3} - \varphi \leq \vartheta_2 \quad \Leftrightarrow \quad \varphi \oplus \vartheta_1 > \varphi \oplus \vartheta_2 \ \
 \text{and} \ \ \vartheta_1 < \vartheta_2.
\]
Consequently, if $\vartheta_1 < \vartheta_2$ and $f_\varphi(\vartheta_1) < f_\varphi(\vartheta_2)$, then $\vartheta_1 < \vartheta_2 < \tfrac{\pi}{3} - \varphi$ or
$\tfrac{\pi}{3} - \varphi \leq \vartheta_1 < \vartheta_2$, and therefore, $f_\varphi$ maps $[\vartheta_1,\vartheta_2]$ onto $[f_\varphi(\vartheta_1),f_\varphi(\vartheta_2)]$ as an increasing bijection.

Setting $\fvk = \Phi(\varphi) \in \nI$, we have, see Figure~\ref{funFI},
\[
F_\fvk (\xi) = \Phi\bigl(f_\varphi\bigl(\Phi^{-1}(\xi)\bigr)\bigr), \qquad \xi \in \nI.
\]
Assume that $\lep < \rep$ and $F_\fvk(\lep) < F_\fvk(\rep)$. Then,
\[
\Phi^{-1}(\lep) < \Phi^{-1}(\rep) \quad  \text{and} \quad  f_\varphi\bigl(\Phi^{-1}(\lep)\bigr) < f_\varphi\bigl(\Phi^{-1}(\rep)\bigr).
\]
Therefore, $f_\varphi$ maps $\bigl[\Phi^{-1}(\lep), \Phi^{-1}(\rep)\bigr]$ onto $\bigl[f_\varphi\bigl(\Phi^{-1}(\lep)\bigr), f_\varphi\bigl(\Phi^{-1}(\rep)\bigr)\bigr]$ as an increasing bijection. Now, $F_\fvk$ restricted to $[\lep, \rep]$ is a composition of three increasing bijections. Thus, (\ref{padxib}) holds.
\end{proof}

\begin{figure}[ht]

\setlength{\abovecaptionskip}{6pt}%
\setlength{\belowcaptionskip}{-8pt}%

\begin{minipage}{0.49\textwidth}

\psfrag{phh}[][]{\begin{picture}(0,0)
            \put(-4,0){\makebox(0,0)[l]{$\tfrac{\varphi}{2}$}}
                        \end{picture}}

\psfrag{php}[][]{\begin{picture}(0,0)
     \put(-13,0){\makebox(0,0)[l]{$\tfrac{\varphi}{2}\!\oplus\!\tfrac{\pi}{6}$}}
                        \end{picture}}

\psfrag{ph}[][]{\begin{picture}(0,0)
     \put(-4,0){\makebox(0,0)[l]{$\varphi$}}
                        \end{picture}}

\psfrag{zx}[][]{\begin{picture}(0,0)
            \put(2,2){\makebox(0,0)[l]{$0$}}
                        \end{picture}}

\psfrag{px3}[][]{\begin{picture}(0,0)
            \put(-4,0){\makebox(0,0)[l]{$\tfrac{\pi}{3}$}}
                        \end{picture}}
\psfrag{px6}[][]{\begin{picture}(0,0)
            \put(-4,0){\makebox(0,0)[l]{$\tfrac{\pi}{6}$}}
                        \end{picture}}
\psfrag{py3}[][]{\begin{picture}(0,0)
            \put(0,-1){\makebox(0,0)[l]{$\tfrac{\pi}{3}$}}
                        \end{picture}}
\psfrag{py6}[][]{\begin{picture}(0,0)
            \put(0,0){\makebox(0,0)[l]{$\tfrac{\pi}{6}$}}
                        \end{picture}}
\resizebox{!}{!}{%
  \includegraphics{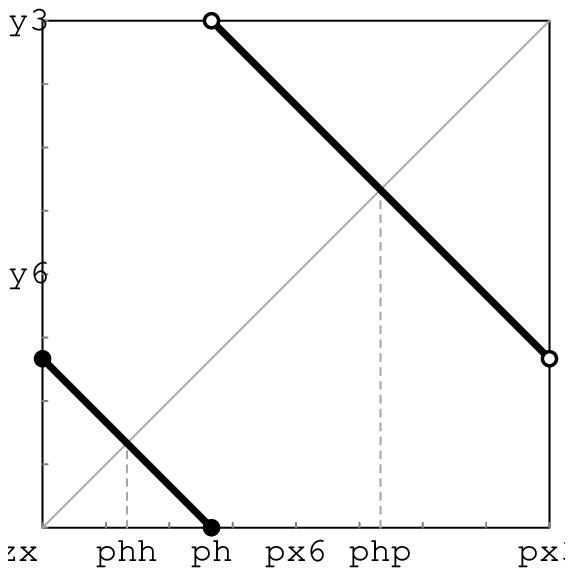}}
    \caption{$g_\varphi$ with $\varphi = \pi/9$}
\label{funF}
\end{minipage} \hfill
\begin{minipage}{0.49\textwidth}

\psfrag{as}[][]{\begin{picture}(0,0)
            \put(-7,1){\makebox(0,0)[l]{$\sqrt[\Bo]{\fvk}$}}
                        \end{picture}}

\psfrag{al}[][]{\begin{picture}(0,0)
            \put(-4,-1){\makebox(0,0)[l]{$\fvk$}}
                        \end{picture}}

\psfrag{asp}[][]{\begin{picture}(0,0)
     \put(-13,1){\makebox(0,0)[l]{$\tfrac{1}{2}\!\BoI\!\sqrt[\Bo]{\fvk}$}}
                        \end{picture}}

\psfrag{zx1}[][]{\begin{picture}(0,0)
            \put(2,2){\makebox(0,0)[l]{$0$}}
                        \end{picture}}

\psfrag{x1}[][]{\begin{picture}(0,0)
            \put(-4,0){\makebox(0,0)[l]{$1$}}
                        \end{picture}}
\psfrag{xh}[][]{\begin{picture}(0,0)
            \put(-4,-1){\makebox(0,0)[l]{$\tfrac{1}{2}$}}
                        \end{picture}}
\psfrag{y1}[][]{\begin{picture}(0,0)
            \put(0,0){\makebox(0,0)[l]{$1$}}
                        \end{picture}}
\psfrag{yh}[][]{\begin{picture}(0,0)
            \put(0,0){\makebox(0,0)[l]{$\tfrac{1}{2}$}}
                        \end{picture}}
\resizebox{!}{!}{%
  \includegraphics{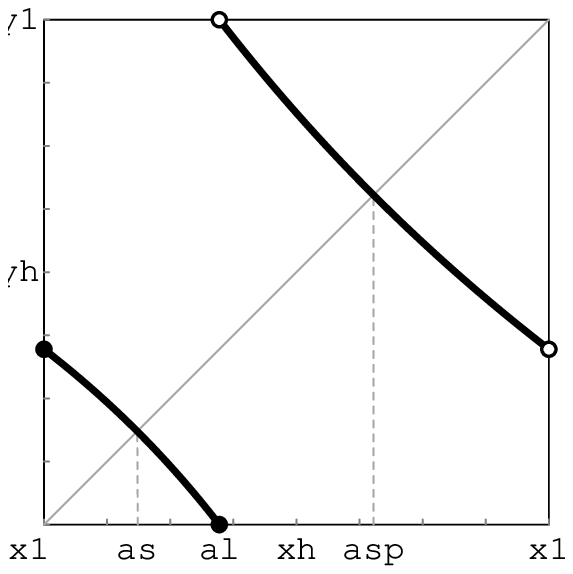}}
    \caption{The function $G_\fvk$}
\label{funG}
\end{minipage}
\end{figure}

\begin{proposition} \label{pimath}
Let \(\fvk, \lep, \rep \in \nI\) and \(\lep < \rep\). Consider the function
\[
G_\fvk(\xi) := \fvk\BoI\TonI{\xi}, \ \xi\in \nI.
\]
\begin{enumerate}[{\rm (a)}]
\item  \label{pimathia}
The function $G_\fvk$ is an involution on $\nI$.
\item  \label{pimathib}
If $G_\fvk(\lep) > G_\fvk(\rep)$, then $G_\fvk$ maps \([\lep, \rep]\) onto \([G_\fvk(\rep), G_\fvk(\lep)]\) as a decreasing bijection.
\item  \label{pimathic}
The fixed points of $G_\fvk$ are $\sqrt[\Bo]{\fvk}$ and $\tfrac{1}{2}\BoI \sqrt[\Bo]{\fvk} = \tfrac{1}{2}\Bo \sqrt[\Bo]{\fvk}$.
\item  \label{pimathid}
The interior of the interval
$\bigl[\sqrt[\Bo]{\fvk}, \tfrac{1}{2}\BoI \sqrt[\Bo]{\fvk}\bigr]$  is mapped onto the exterior of this interval in $\nI$.
\end{enumerate}
\end{proposition}
\begin{proof}
The statement (\ref{pimathia}) is clear. To prove (\ref{pimathib})
let $\varphi \in \bigl(0,\tfrac{\pi}{3}\bigr)$ and consider the function $g_\varphi(\vartheta) = \varphi \ominus \vartheta$, see Figure~\ref{funF}. Since
\[
g_\varphi(\vartheta) = \left\{ \begin{array}{ll}
 \varphi - \vartheta & \
    \text{if} \quad 0 \leq \vartheta \leq \varphi, \\[4pt]
 \varphi - \vartheta + \tfrac{\pi}{3} & \  \text{if} \quad \varphi < \vartheta < \tfrac{\pi}{3},
 \end{array}
\right.
\]
we have the equivalence
\[
\vartheta_1 \leq \varphi < \vartheta_2 \quad \Leftrightarrow \quad \varphi \ominus \vartheta_1 < \varphi \ominus \vartheta_2 \ \
 \text{and} \ \ \vartheta_1 < \vartheta_2.
\]
Consequently, if $\vartheta_1 < \vartheta_2$ and $g_\varphi(\vartheta_1) > g_\varphi(\vartheta_2)$, then $\vartheta_1 < \vartheta_2 \leq \varphi$ or
$\varphi < \vartheta_1 < \vartheta_2$, and therefore, $g_\varphi$ maps $[\vartheta_1,\vartheta_2]$ onto $[g_\varphi(\vartheta_2),g_\varphi(\vartheta_1)]$ as a decreasing bijection.

Setting $\fvk = \Phi(\varphi)$ we have, see Figure~\ref{funG},
\[
G_\fvk(\xi) = \Phi\bigl(g_\varphi\bigl(\Phi^{-1}(\xi)\bigr)\bigr), \qquad \xi \in \nI.
\]
Assume that $\lep < \rep$ and $G_\fvk(\lep) > G_\fvk(\rep)$. Then
\[
\Phi^{-1}(\lep) < \Phi^{-1}(\rep) \quad  \text{and} \quad  g_\varphi\bigl(\Phi^{-1}(\lep)\bigr) > g_\varphi\bigl(\Phi^{-1}(\rep)\bigr).
\]
Therefore $g_\varphi$ maps $\bigl[\Phi^{-1}(\lep), \Phi^{-1}(\rep)\bigr]$ onto $\bigl[g_\varphi\bigl(\Phi^{-1}(\rep)\bigr), g_\varphi\bigl(\Phi^{-1}(\lep)\bigr)\bigr]$ as a decreasing bijection. Now, $G_\fvk$ restricted to $[\lep, \rep]$ is a composition of two increasing bijections and one decreasing bijection. Thus (\ref{pimathib}) holds.

The statement (\ref{pimathic}) was proved at the beginning of this section. To prove (\ref{pimathid}) we use the function $g_\varphi$ again with $\varphi = \Phi^{-1}(\fvk)$. The fixed points of $g_\varphi$ are $\tfrac{\varphi}{2}$ and $\tfrac{\varphi}{2}\oplus\tfrac{\pi}{6}$, see Figure~\ref{funF}. It is clear that $g_\varphi$ maps $\bigl(\tfrac{\varphi}{2},\varphi\bigr]$ to $\bigl[0,\tfrac{\varphi}{2}\bigr)$. Also, $g_\varphi$ maps $\bigl(\varphi, \tfrac{\varphi}{2}+\tfrac{\pi}{6}\bigr)$ to   $\bigl(\tfrac{\varphi}{2}+\tfrac{\pi}{6},1\bigr)$. That is, $g_\varphi$ maps the interior of $\bigl[\tfrac{\varphi}{2},\tfrac{\varphi}{2}+\tfrac{\pi}{6}\bigr]$ onto its exterior. The statement (\ref{pimathid}) now follows from the the fact that $\Phi^{-1}$ maps $\bigl[\sqrt[\Bo]{\fvk}, \tfrac{1}{2}\BoI \sqrt[\Bo]{\fvk}\bigr]$ onto $\bigl[\tfrac{\varphi}{2},\tfrac{\varphi}{2}+\tfrac{\pi}{6}\bigr]$ and
$\Phi$ maps the exterior of  $\bigl[\tfrac{\varphi}{2},\tfrac{\varphi}{2}+\tfrac{\pi}{6}\bigr]$ onto the exterior of  $\bigl[\sqrt[\Bo]{\fvk}, \tfrac{1}{2}\BoI \sqrt[\Bo]{\fvk}\bigr]$.
\end{proof}

We end this section with a few formulas connecting the $\Bo$-operation and the function $\pn{\,\cdot\,}$. The following three identities for $\pn{\,\cdot\,}$ are verified by simplification of the left hand sides:
\begin{align} \label{eqrt}
\llparenthesis \rho\Bo \tau \rrparenthesis & = \frac{1}{|1-\rho\tau|} \llparenthesis \rho \rrparenthesis \, \llparenthesis \tau\rrparenthesis, \quad \rho, \tau \in \nR, \ \rho\tau\neq 1, \\ \label{eqrinfty}
\llparenthesis \rho\Bo \infty \rrparenthesis & = \frac{1}{|\rho|} \llparenthesis \rho\rrparenthesis, \quad \rho \in \nR\setminus\{0\},\\
\label{eqrinv}
\llparenthesis\Ton{\rho}\rrparenthesis & = \frac{1}{|1-\rho|} \llparenthesis \rho\rrparenthesis, \quad \rho \in \nR\setminus\{1\}.
\end{align}
From (\ref{eqrt}), (\ref{eqrinfty}) and (\ref{eqrinv}) we obtain that, whenever the right hand sides are defined, the following identities hold as well:
{\allowdisplaybreaks
\begin{align} \nonumber
\llparenthesis \rho\Bo{\tau}\Bo \zeta\rrparenthesis & = \frac{1}{|1-\rho\tau-\tau\zeta-\zeta\rho+\rho\tau\zeta|} \,\llparenthesis \rho \rrparenthesis \, \llparenthesis \tau\rrparenthesis \, \llparenthesis \zeta\rrparenthesis , \quad \rho, \tau, \zeta \in \nR, \\
\label{eqthreeb}
\llparenthesis \rho\Bo\Ton{\tau}\Bo \zeta\rrparenthesis & = \frac{1}{|1-\tau+\rho\tau+\tau\zeta-\zeta\rho|} \,\llparenthesis \rho \rrparenthesis \, \llparenthesis \tau\rrparenthesis \, \llparenthesis \zeta\rrparenthesis , \quad \rho, \tau, \zeta \in \nR,  \\
\label{eqtwoinfty}
\llparenthesis \rho\Bo\Ton{\tau}\Bo \infty \rrparenthesis
  & = \frac{1}{|\tau-\rho|} \,\llparenthesis \rho \rrparenthesis \, \llparenthesis \tau\rrparenthesis , \quad \rho, \tau \in \nR, \\
\label{eqtwoone}
\llparenthesis \rho\Bo 1 \Bo \zeta \rrparenthesis
  & = \frac{1}{|\rho+\zeta-1|} \,\llparenthesis \rho \rrparenthesis \, \llparenthesis \zeta \rrparenthesis , \quad \rho, \tau, \zeta \in \nR.
\end{align}}

\section{Groups and reflections}

Let $\rho \in \nRc$ and $\vartheta_\rho = \Phi^{-1}(\rho)$. Denote by $\Ref(\vartheta_\rho)$ the matrix with respect to \(\{\Vb{p}_0, \Vb{q}_0\}\) of the reflection across the line determined by the vector $\Vb{p}_\rho$ in the plane spanned by the vectors \(\Vb{p}_0, \Vb{q}_0.\) Then, the matrix with respect to the basis \(\{\Vb{p}_0, \Vb{q}_0, \Vb{r}\}\) of the reflection induced by $\Vs{M}_\rho$ is
\begin{equation} \label{eqrefl}
\left[\!\! \begin{array}{cc}
\Ref(\vartheta_\rho) & \!\!\!\!\! \begin{array}{c}
0 \\ 0
\end{array} \\
\begin{array}{cc}
\! 0 & \, 0
\end{array} &\!\!\!\!\! 1
\end{array}\!\! \right] = \Vs{Q}^\top \Vs{M}_\rho \Vs{Q},
\end{equation}
where $\Vs{Q}$ is the orthogonal matrix whose columns are the vectors $\Vb{p}_0, \Vb{q}_0, \Vb{r}$. From now on, the $3\times 3$ matrix in \eqref{eqrefl} will be identified with its top left corner $\Ref(\vartheta_\rho)$. In the same spirit, we denote by  $\Rot(\vartheta)$ the matrix with respect to the basis \(\{\Vb{p}_0, \Vb{q}_0, \Vb{r}\}\) of the  counterclockwise rotation about the vector $\Vb{r}$ by the angle $\vartheta$.  Familiar formulas connecting coordinate rotations and reflections in the plane  extend to this setting:
\begin{equation} \label{eqrr}
\Ref(\vartheta) \, \Ref(\varphi) = \Rot(2(\vartheta - \varphi)) \quad \text{and} \quad \Rot(\vartheta) \, \Ref(\varphi) = \Ref(\varphi + \vartheta/2).
 \end{equation}

\begin{theorem} \label{tcomp}
Let $\rho,  \tau, \zeta \in \nRc$. We have the following matrix identity:
 $$
\Vs{M}_\zeta \Vs{M}_\tau \Vs{M}_\rho = \Vs{M}_{\rho\Bo \Ton{\tau} \Bo \zeta}.
 $$
\end{theorem}
\begin{proof}
Let $\rho, \tau,\zeta \in \nRc$ and $\vartheta_\rho = \Phi^{-1}(\rho), \vartheta_\tau = \Phi^{-1}(\tau), \vartheta_\zeta = \Phi^{-1}(\zeta)$. Since by definition of $\Bo$ the mapping $\Phi$ is an isomorphism between the groups $\Bigl(\bigl[-\frac{\pi}{3},\frac{2\pi}{3}\bigr),\oplus\Bigl)$ and $(\nRc,\Bo)$ we have that $\Phi(\vartheta_\rho - \vartheta_\tau + \vartheta_\zeta) = \rho\Bo \Ton{\tau}\Bo \zeta$. Together with \eqref{eqrefl} and \eqref{eqrr} this yields the following equalities
\begin{align*}
\Vs{M}_\rho \Vs{M}_\tau \Vs{M}_\zeta & = \Vs{Q} \Ref(\vartheta_\rho) \, \Ref(\vartheta_\tau) \, \Ref(\vartheta_\zeta) \, \Vs{Q}^\top \\
 & = \Vs{Q} \Rot(2(\vartheta_\rho - \vartheta_\tau)) \Ref(\vartheta_\zeta) \,\Vs{Q}^\top \\
 & = \Vs{Q} \Ref(\vartheta_\rho-\vartheta_\tau+\vartheta_\zeta)\,\Vs{Q}^\top \\
 & = \Vs{Q} \Ref(\vartheta_{\rho\Bo \Ton{\tau} \Bo \zeta})\,\Vs{Q}^\top \\
  & = \Vs{M}_{\rho\Bo \Ton{\tau} \Bo \zeta}. \qedhere
\end{align*}
\end{proof}

Recall that the cone $\cQ$ was introduced in Section~\ref{scone} and the vector $\Vb{r}$ was introduced in Section~\ref{refm}.
\begin{proposition} \label{prsu}
Let $\rho, \sigma \in \nRc$ and let $\Vb{t}\in \cQ$ be such that $\Vb{t}$ and $\Vb{r}$ are linearly independent. Then $\Vs{M}_\sigma \Vb{t} = \lambda  \Vs{M}_\rho \Vb{t}$ if and only if $\lambda = 1$ and $\sigma = \rho$.
\end{proposition}
\begin{proof}
Since both vectors $\Vs{M}_\sigma \Vb{t}$ and $\Vs{M}_\rho \Vb{t}$ are in $\cQ$, $\lambda$ must be positive. Applying $\Vs{M}_{\rho}$ to both sides of $\Vs{M}_\sigma \Vb{t} = \lambda  \Vs{M}_{\rho} \Vb{t}$ we get $\Vs{M}_\rho \Vs{M}_\sigma \Vb{t} = \lambda \Vb{t}$. By \eqref{eqrefl} and \eqref{eqrr} $\Vs{M}_\rho \Vs{M}_\sigma = \Vs{Q} \Rot (2(\vartheta_{\rho} - \vartheta_{\sigma}))\Vs{Q}^\top$ and since only the identity rotation has the positive eigenvalue $1$, we conclude that $\vartheta_\rho = \vartheta_\sigma$ and $\lambda = 1$.
\end{proof}

Let $\Vb{t} =[a' \ b' \ c']^\top \in \nR^3$ be such that $\Vb{t}$ and $\Vb{r}$ are linearly independent. Set
\begin{equation*}
\varkappa_{\Vb{t}} := \left\{ \begin{array}{lc}
\dfrac{a'-b'}{a'-c'} &  \text{if} \quad a'\neq c', \\[10pt]
\infty & \text{if} \quad a'=c'.
\end{array} \right.
\end{equation*}

\begin{theorem} \label{pgamma}
Let $\Vb{t}$ and $\Vb{v}$ be nonzero vectors in $\nR^3$, neither of which is a multiple of $\Vb{r}$. Then  $\Vs{M}_\rho \Vb{t} = \lambda \Vb{v}$ has a unique solution for $\lambda \in \nR$ and $\rho \in \nRc$ if and only if $\|\Vb{v}\| (\Vb{t}\cdot\Vb{r}) = \|\Vb{t}\| (\Vb{v}\cdot\Vb{r})$. In this case,
\begin{equation*}
\lambda = \dfrac{\|\Vb{t}\|}{\|\Vb{v}\|}, \qquad \text{and}  \qquad \rho =
\sqrt[\Bo]{\varkappa_{\Vb{v}}} \Bo  \sqrt[\Bo]{\varkappa_{\Vb{t}}} \quad  \text{or} \quad  \rho =
\sqrt[\Bo]{\varkappa_{\Vb{v}}} \Bo \sqrt[\Bo]{\varkappa_{\Vb{t}}} \Bo 2.
\end{equation*}
In particular, if $\Vb{v}=\Vb{t}$, then $\rho = \varkappa_{\Vb{t}}$.
\end{theorem}

\begin{proof}
It is a lengthy but straightforward calculation to verify that $\Vs{M}_{\varkappa_{\Vb{t}}} \Vb{t} = \Vb{t}$. The uniqueness follows from Proposition~\ref{prsu}. So we only need to prove that $\Vs{M}_\rho \Vb{t} = \lambda \Vb{v}$ is equivalent to $\|\Vb{v}\| (\Vb{t}\cdot\Vb{r}) = \|\Vb{t}\| (\Vb{v}\cdot\Vb{r})$. For simplicity, and without loss of generality, we assume that $\Vb{t}$ and $\Vb{v}$ are unit vectors. Since $\Vs{M}_\rho$ is a reflection, it preserves length, and therefore $\lambda=1$. Then, $\Vs{M}_\rho \Vb{t} = \Vb{v}$ implies $\,\Vb{t}\cdot\Vb{r} = \Vb{v}\cdot\Vb{r}$, since $\Vs{M}_\rho$ is a unitary mapping and $\Vs{M}_\rho\Vb{r} = \Vb{r}$. To prove the converse, assume $\,\Vb{t}\cdot\Vb{r} = \Vb{v}\cdot\Vb{r}$. By definition, we have
\begin{equation*}
\sqrt[\Bo]{\varkappa_{\Vb{t}}} \Bo \sqrt[\Bo]{\varkappa_{\Vb{t}}} = \varkappa_{\Vb{t}} \qquad \text{and} \qquad \sqrt[\Bo]{\varkappa_{\Vb{v}}} \Bo \sqrt[\Bo]{\varkappa_{\Vb{v}}} = \varkappa_{\Vb{v}}.
\end{equation*}
By Theorem~\ref{tcomp}, the last two equations imply $\Vs{M}_{\!\sqrt[\Bo]{\varkappa_{\Vb{t}}}}\Vs{M}_0 \Vs{M}_{\!\sqrt[\Bo]{\varkappa_{\Vb{t}}}} = \Vs{M}_{\!\varkappa_{\Vb{t}}}$ and $\Vs{M}_{\!\sqrt[\Bo]{\varkappa_{\Vb{v}}}}\Vs{M}_0 \Vs{M}_{\!\sqrt[\Bo]{\varkappa_{\Vb{t}}}} = \Vs{M}_{\!\varkappa_{\Vb{v}}}$. Since \(\Vs{M}_{\!\varkappa_{\Vb{t}}} \Vb{t} = \Vb{t}\) and \(\Vs{M}_{\!\varkappa_{\Vb{v}}} \Vb{v} = \Vb{v}\), we have
\begin{equation} \label{eqM0}
 \Vs{M}_0 \Vs{M}_{\!\raisebox{-2pt}{\scriptsize$\sqrt[\Bo]{\varkappa_{\Vb{t}}}$}} \Vb{t}
  = \Vs{M}_{\!\raisebox{-2pt}{\scriptsize$\sqrt[\Bo]{\varkappa_{\Vb{t}}}$}} \Vb{t}
   \qquad \text{and} \qquad
   \Vs{M}_0 \Vs{M}_{\!\raisebox{-2pt}{\scriptsize$\sqrt[\Bo]{\varkappa_{\Vb{v}}}$}} \Vb{v}
    = \Vs{M}_{\!\raisebox{-2pt}{\scriptsize$\sqrt[\Bo]{\varkappa_{\Vb{v}}}$}} \Vb{v}.
\end{equation}
This implies that the unit vectors \(\Vs{M}_{\!\raisebox{-2pt}{\scriptsize$\sqrt[\Bo]{\varkappa_{\Vb{t}}}$}} \Vb{t}\) and \(\Vs{M}_{\!\raisebox{-2pt}{\scriptsize$\sqrt[\Bo]{\varkappa_{\Vb{v}}}$}} \Vb{v}\) are in the plane which is invariant under $\Vs{M}_0$. The assumption $\,\Vb{t}\cdot\Vb{r} = \Vb{v}\cdot\Vb{r}$ yields $\,\bigl( \Vs{M}_{\!\raisebox{-2pt}{\scriptsize$\sqrt[\Bo]{\varkappa_{\Vb{t}}}$}} \Vb{t} \bigr)\cdot\Vb{r} = \bigl(\Vs{M}_{\!\raisebox{-2pt}{\scriptsize$\sqrt[\Bo]{\varkappa_{\Vb{v}}}$}} \Vb{v} \bigr)\cdot\Vb{r}$. Therefore, either
\[
 \Vs{M}_{\!\raisebox{-2pt}{\scriptsize$\sqrt[\Bo]{\varkappa_{\Vb{v}}}$}} \Vb{v}
  = \Vs{M}_{\!\raisebox{-2pt}{\scriptsize$\sqrt[\Bo]{\varkappa_{\Vb{t}}}$}} \Vb{t}
   \qquad \text{or} \qquad
 \Vs{M}_{\!\raisebox{-2pt}{\scriptsize$\sqrt[\Bo]{\varkappa_{\Vb{v}}}$}} \Vb{v} = \Vs{M}_2 \Vs{M}_{\!\raisebox{-2pt}{\scriptsize$\sqrt[\Bo]{\varkappa_{\Vb{t}}}$}} \Vb{t}.
\]
In the first case, we substitute the equality in \eqref{eqM0} and apply  \(\Vs{M}_{\!\sqrt[\Bo]{\varkappa_{\Vb{t}}}}\), while in the second case we just apply  \(\Vs{M}_{\!\sqrt[\Bo]{\varkappa_{\Vb{v}}}}\) to get
\[
 \Vs{M}_{\!\raisebox{-2pt}{\scriptsize$\sqrt[\Bo]{\varkappa_{\Vb{v}}}$}} \Vs{M}_0
  \Vs{M}_{\!\raisebox{-2pt}{\scriptsize$\sqrt[\Bo]{\varkappa_{\Vb{t}}}$}} \Vb{t} = \Vb{v}
 \qquad \text{or} \qquad
  \Vs{M}_{\!\raisebox{-2pt}{\scriptsize$\sqrt[\Bo]{\varkappa_{\Vb{v}}}$}} \Vs{M}_2 \Vs{M}_{\!\raisebox{-2pt}{\scriptsize$\sqrt[\Bo]{\varkappa_{\Vb{t}}}$}} \Vb{t} = \Vb{v}.
\]
Now Theorem~\ref{tcomp} and Proposition~\ref{prsu} yield the claim.
\end{proof}

\section{Similarity of triangles in $\CT_{\mathbb R} (T)$} \label{sprop}

For the remainder of this paper, we will use the following notation. Given an oriented triangle $T = (a,b,c)$, $\Vb{t}$ will denote the vector in $\nR^3$ whose components are the squares of the sides of $T$, that is $\Vb{t}=[a^2 \ b^2 \ c^2]^{\top}$. We also set $\varkappa_{_T} : = \varkappa_{\Vb{t}}$. The following relationship, which follows from \eqref{eqsqMsq}, is important in the reasoning below and we will use it without explicitly mentioning it:
\[
 V = \CT_\rho(T) \quad \Leftrightarrow \quad  \Vb{v} = \llparenthesis \rho \rrparenthesis^2 \Vs{M}_\rho \Vb{t} .
\]

In the rest of the paper we use the algebraic set-up established in the previous four sections to investigate the structure of the family of Ceva's triangles $\CT_{\mathbb R} (T)$. First, we establish a simple relationship between Ceva's triangles of reversely congruent oriented triangles.

\begin{proposition} \label{pinvcon}
Let $T$ and $V$ be reversely congruent oriented triangles. Then, for every $\rho$ in $\nRc$, the triangles $\CT_\rho (T)$ and $\CT_{1-\rho} (V)$ are reversely congruent.
\end{proposition}
\begin{proof}
Since $T$ and $V$ are reversely congruent, there exists $\sigma \in \nS$ such that $\Vb{v} = \Vs{M}_\sigma \Vb{t}$. Since $\sigma = \rho\Bo\Ton{(1\Bo\Ton{\sigma})}\Bo 1\Bo \Ton{\rho}$, we have $\Vb{v} = \Vs{M}_{1\Bo \Ton{\rho}} \Vs{M}_{1\Bo\Ton{\sigma}} \Vs{M}_{\rho}\Vb{t}$, and consequently $\Vs{M}_{1\Bo \Ton{\rho}} \Vb{v} = \Vs{M}_{1\Bo\Ton{\sigma}} \Vs{M}_{\rho}\Vb{t}$. Since $\nS$ is a subgroup of $(\nR, \Bo)$, $1\Bo\Ton{\sigma} \in \nS$. Since $1\Bo \Ton{\rho} = 1- \rho$, the triangles $\CT_\rho (T)$ and $\CT_{1-\rho} (V)$ are reversely similar. Since
$\pn{1-\rho} =\pn{\rho}$, the ratio of similarity is $1$.
\end{proof}

\begin{proposition} \label{prepeat}
Let $T$ be an oriented triangle and $\rho \in \nR$. Then the triangles
$\CT_\rho \bigl(\CT_\rho(T)\bigr)$ and $T$ are directly similar in the ratio $\pn{\rho}^2$. The triangle $\CT_\infty \bigl(\CT_\infty(T)\bigr)$ is directly congruent to $T$. More precisely, $\CT_\rho^2(T) = \CT_\rho \bigl(\CT_\rho(T)\bigr) = \pn{\rho}^2 T$ and $\CT_\infty \bigl(\CT_\infty(T)\bigr) = T$.
\end{proposition}
\begin{proof}
The squares of the sides of $\CT_\rho(T)$ are the components of the vector $\pn{\rho}^4 \Vs{M}_\rho \Vb{t}$. Hence, the squares of the sides of $\CT_\rho \bigl(\CT_\rho(T)\bigr)$ are the components of the vector $\pn{\rho}^4 \Vs{M}_\rho \Vs{M}_\rho \Vb{t}=\pn{\rho}^4\Vb{t}$. This implies the first statement. The second statement is straightforward application of the definition.
\end{proof}

Letting $\rho=1/2$ in the preceding proposition, we recover a fact mentioned in the Introduction: the median of the median triangle is similar to the original triangle in the ratio $3/4$. Notice also that combining Propositions~\ref{pinvcon} and~\ref{prepeat} proves the converse of Proposition~\ref{pinvcon}. Another immediate consequence of Propositions~\ref{pinvcon} and~\ref{prepeat} is the following corollary.

\begin{corollary}\label{rho-T-V}
Let  $\rho\in \mathbb R$. The oriented triangles $T$ and $V$ are directly similar with the ratio of similarity $l$ if and only if $\CT_\rho (T)$ and $\CT_\rho (V)$ are directly similar with the ratio of similarity $l$. In particular, for $l > 0$, $\CT_\rho (l\, T) =  l \, \CT_\rho (T)$ and $\CT_\rho^{2n}(T) = \pn{\rho}^{2n} T$ for all $n \in \nN$.
\end{corollary}

It is worth noting that, given an oriented triangle $T$, the study of the
similarity properties of the family $\CT_{\mathbb R} (T)$ can be
reduced to that of the subfamily  $\CT_{\mathbb I} (T)$. Recall that
${\mathbb I} = [0,1)$ and $\nS = \{0,1,\infty\}$.

\begin{theorem}[Reduction to $\CT_{\mathbb I} (T)$]
  \label{reduction 01}
Let  $\rho, \tau \in \nRc$, and let $T$ be a non-equilateral oriented triangle. The triangles $\CT_{\rho}(T)$ and $\CT_{\tau}(T)$ are directly similar if and only if $\tau \in \rho \Bo \nS$. For $\rho\in\nR\setminus\{0\}$, the ratio of similarity of $\CT_{\rho}(T)$ to $\CT_{\rho\Bo\infty}(T)$ is $|\rho|$. For $\rho\in\nR\setminus\{1\}$, the ratio of similarity of $\CT_{\rho}(T)$ to $\CT_{\rho\Bo 1}(T)$ is $|1-\rho|$.
\end{theorem}
\begin{proof}
The triangles $\CT_{\rho}(T)$, $\CT_{\tau}(T)$ are directly similar if and only if there exists $\phi \in \{0,2\pi/3,-2\pi/3\}$ such that the counterclockwise  rotation of $\Vs{M}_\rho \Vb{t}$ about $\Vb{r}$ by $\phi$ coincides with $\Vs{M}_\tau \Vb{t}$. The last condition restated in matrix notation reads
\[
\Vs{M}_\tau \Vb{t} = \Vs{Q} \Rot(\phi) \Vs{Q}^\top \Vs{M}_\rho \Vb{t}
\]
Next, we use the following identities:
\begin{align*}
\Vs{M}_0 \Vs{M}_1 & =\Vs{Q} \Ref(0)\Ref(\pi/3)\,\Vs{Q}^\top = \Vs{Q}\Rot(2 \pi/3)\,\Vs{Q}^\top, \\
\Vs{M}_1 \Vs{M}_0 & =\Vs{Q} \Ref(\pi/3)\Ref(0)\,\Vs{Q}^\top = \Vs{Q}\Rot(-2 \pi/3)\,\Vs{Q}^\top.
\end{align*}
If $\phi = 0$, then $\tau = \rho$ by Proposition~\ref{prsu}. If $\phi = 2 \pi/3$, then by Theorem~\ref{tcomp} and Proposition~\ref{prsu} we have $\tau =  0\Bo\rho \Bo \Ton{1} = \rho \Bo \infty$. Similarly, if $\phi = -2 \pi/3$, then $\tau =  1\Bo\rho\Bo\Ton{0} = \rho\Bo 1$. To calculate the ratios of similarity, consider $\tau = \rho\Bo \infty$ and $\rho \neq 0$  first. The squares of the sides of $\CT_\tau(T)$ are the components of the vector $\llparenthesis \rho\Bo\infty \rrparenthesis^2 \Vs{M}_{\rho\Bo \infty} \Vb{t}$, while the squares of the sides of $\CT_\rho(T)$ are the components of the vector $\llparenthesis \rho \rrparenthesis^2 \Vs{M}_\rho \Vb{t}$. Since, by \eqref{eqrinfty}, $\llparenthesis \rho\Bo\infty \rrparenthesis = \llparenthesis \rho\rrparenthesis/|\rho|$, the ratio of similarity of  $\CT_{\rho}(T)$ to $\CT_{\rho\Bo 1}(T)$ is $|\rho|$. The remaining claim is proved using \eqref{eqrt}.
\end{proof}

The following result is an extension of Parry's characterization of automedian triangles, see \cite{Parry}.

\begin{corollary} \label{comega}
Let $T$ be a non-equilateral oriented triangle. Then $T$ is directly similar to $\CT_{\rho}(T)$ if and only if $\rho \in \varkappa_{_{T}} \Bo \nS$. The ratio of similarity of  $\CT_{\rho}(T)$ to $T$ is $\pn{\rho}$ with $\rho \in \varkappa_{_T}\Bo\nS$. In particular, $T$ is automedian if and only if $1/2 \in \varkappa_{_T}\Bo\nS$.
\end{corollary}
\begin{proof}
By Theorem~\ref{pgamma}, $\Vs{M}_{\varkappa_{\Vb{t}}} \Vb{t} = \Vb{t}$. Therefore, $T$ and $\CT_{\varkappa_{_T}}(T)$ are directly similar. By Theorem~\ref{reduction 01},  $\CT_{\rho}(T)$ is directly similar to $\CT_{\varkappa_{_T}}(T)$ if and only if $\rho \in \varkappa_{_T} \Bo \nS$.
\end{proof}

\begin{figure}[ht]

\setlength{\abovecaptionskip}{6pt}%
\setlength{\belowcaptionskip}{-8pt}%

\psfrag{A}[][]{\begin{picture}(0,0)
            \put(-18,-2){\makebox(0,0)[l]{$A$}}
                        \end{picture}}
\psfrag{Ap}[][]{\begin{picture}(0,0)
            \put(0,-2){\makebox(0,0)[l]{$A_{\text{{\tiny $1/2$}}}$}}
                        \end{picture}}
\psfrag{Aq}[][]{\begin{picture}(0,0)
            \put(-8,-2){\makebox(0,0)[l]{$A_{\text{{\tiny $-1$}}}$}}
                        \end{picture}}
\psfrag{Ar}[][]{\begin{picture}(0,0)
            \put(-6,-3){\makebox(0,0)[l]{$A_{\text{{\tiny $2$}}}$}}
                        \end{picture}}

\psfrag{B}[][]{\begin{picture}(0,0)
            \put(1,-2){\makebox(0,0)[l]{$B$}}
                        \end{picture}}

\psfrag{Bp}[][]{\begin{picture}(0,0)
            \put(-7,1){\makebox(0,0)[l]{$B_{\text{{\tiny $1/2$}}}$}}
                        \end{picture}}

\psfrag{Bq}[][]{\begin{picture}(0,0)
            \put(-5,-2){\makebox(0,0)[l]{$B_{\text{{\tiny $-1$}}}$}}
                        \end{picture}}
\psfrag{Br}[][]{\begin{picture}(0,0)
            \put(-5,2){\makebox(0,0)[l]{$B_{\text{{\tiny $2$}}}$}}
                        \end{picture}}

\psfrag{C}[][]{\begin{picture}(0,0)
            \put(0,-1){\makebox(0,0)[l]{$C$}}
                        \end{picture}}

\psfrag{Cp}[][]{\begin{picture}(0,0)
            \put(-27,1){\makebox(0,0)[l]{$C_{\text{{\tiny $1/2$}}}$}}
                             \end{picture}}

\psfrag{Cq}[][]{\begin{picture}(0,0)
            \put(-5,2){\makebox(0,0)[l]{$C_{\text{{\tiny $-1$}}}$}}
                             \end{picture}}
\psfrag{Cr}[][]{\begin{picture}(0,0)
            \put(-5,-3){\makebox(0,0)[l]{$C_{\text{{\tiny $2$}}}$}}
                             \end{picture}}

\psfrag{App}[][]{\begin{picture}(0,0)
            \put(-6,-1){\makebox(0,0)[l]{$A_{\text{{\tiny $1/2$}}}'$}}
                        \end{picture}}

\psfrag{Apq}[][]{\begin{picture}(0,0)
            \put(-8,-3){\makebox(0,0)[l]{$A_{\text{{\tiny $-1$}}}'$}}
                        \end{picture}}
\psfrag{Apr}[][]{\begin{picture}(0,0)
            \put(-5,2){\makebox(0,0)[l]{$A_{\text{{\tiny $2$}}}'$}}
                        \end{picture}}

\resizebox{!}{!}{%
  \includegraphics{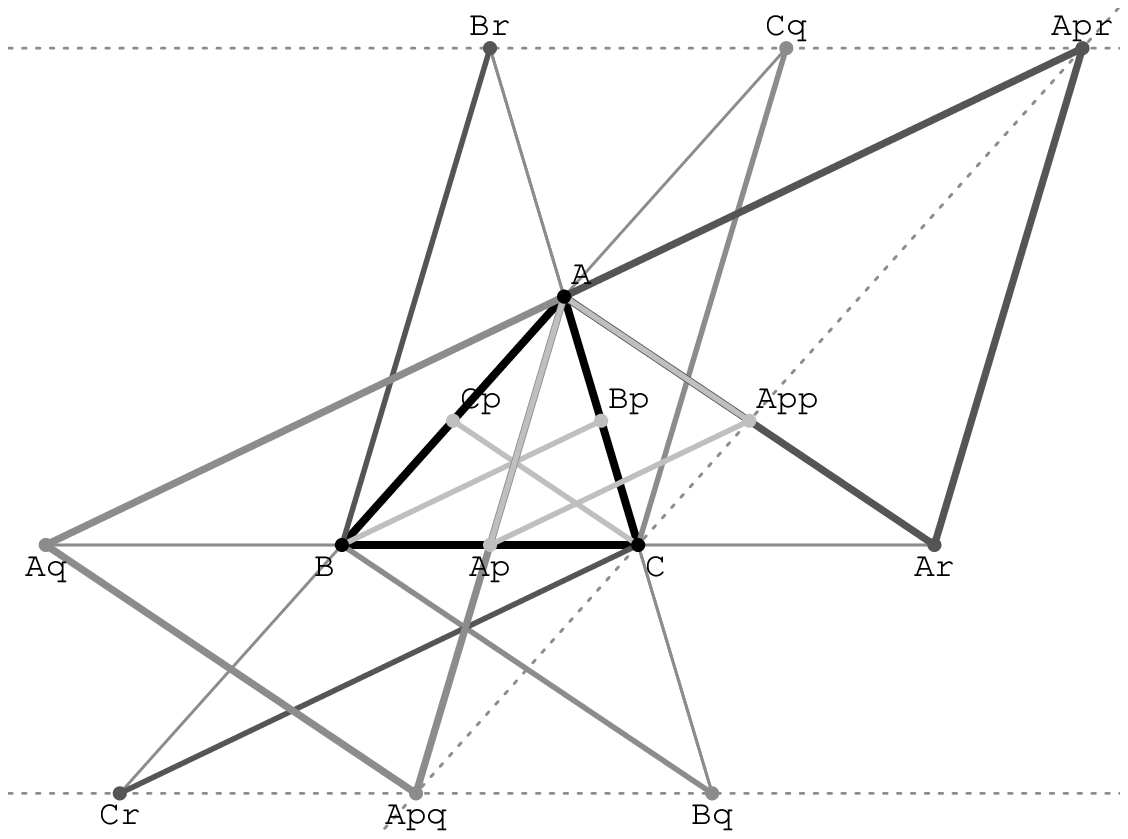}}
   \caption{
   $AA_{\text{{\tiny $1/2$}}}A_{\text{{\tiny $1/2$}}}'$
   is reversely congruent to $\CT_{1/2}(T)$;\\
   $AA_{\text{{\tiny $-1$}}}A_{\text{{\tiny $-1$}}}'$ and $AA_{\text{{\tiny $2$}}}A_{\text{{\tiny $2$}}}'$ are directly congruent and reversely congruent to $\CT_{-1}(T)$ and $\CT_{2}(T)$}
\label{f3med}

\end{figure}

\begin{corollary}
The only two pairs of directly congruent triangles in $\CT_{\nR}(T)$ are   $\CT_0(T), \CT_1(T)$ and $\CT_{-1}(T), \CT_2(T)$.
\end{corollary}
\begin{proof}
By Theorem~\ref{reduction 01}, for $\CT_\rho(T)$ and $\CT_\sigma(T)$ to be directly congruent we must have $\sigma = \rho\Bo\infty$ or $\sigma = \rho\Bo 1$, and $|\rho| = 1$ or $|1-\rho| = 1$, respectively. Clearly, the only candidates are $\rho \in \{-1,0,1,2\}$. The corresponding $\sigma$-s are $2, 1, 0, -1$, respectively. Each of the triangles $\CT_0(T), \CT_1(T), \CT_\infty(T)$ is reversely congruent to $T$. Each of the triangles $\CT_{-1}(T), \CT_2(T)$ is directly similar to the median triangle $\CT_{1/2}(T)$, with the ratio of similarity $2$, see Figure~\ref{f3med}.
\end{proof}

\begin{lemma} \label{lcomp}
Let $\rho,\tau,\zeta \in \nRc$ and let $T$ be an oriented non-equilateral triangle. Then the triangle $\CT_\zeta\bigl(\CT_\tau\bigl(\CT_\rho(T)\bigr)\bigr)$ is directly similar to $\CT_{\rho \Bo\Ton{\tau}\Bo\zeta}(T)$, with the ratio of similarity
\begin{equation} \label{eqrasi}
\frac{\pn{\rho}\pn{\tau}\pn{\zeta}}{\pn{\rho \Bo\Ton{\tau}\Bo\zeta}}.
\end{equation}
\end{lemma}
\begin{proof}
Set $V = \CT_\zeta\bigl(\CT_\tau\bigl(\CT_\rho(T)\bigr)\bigr)$ and $U = \CT_{\rho \Bo\Ton{\tau}\Bo\zeta}(T)$. By Theorem~\ref{tcomp}, we have $\Vs{M}_\zeta \Vs{M}_\tau \Vs{M}_\rho = \Vs{M}_{\rho \Bo\Ton{\tau}\Bo\zeta}$. Therefore, $\Vb{v} = \pn{\zeta}^2 \pn{\tau}^2 \pn{\rho}^2 \Vs{M}_{\rho \Bo\Ton{\tau}\Bo\zeta} \Vb{t}.$ By \eqref{eqsqMsq}, $\Vb{u} = \pn{\rho \Bo\Ton{\tau}\Bo\zeta}^2
\Vs{M}_{\rho \Bo\Ton{\tau}\Bo\zeta} \Vb{t}.$ Thus, the triangle $V$ is directly similar to the triangle $U$ and the ratio of similarity is given by \eqref{eqrasi}.
Assume that $\rho \Bo\Ton{\tau}\Bo\zeta \neq \infty$. Then, if $\rho,\tau,\zeta \in \nR$, the ratio of similarity simplifies to $|1-\tau +\rho\tau+\tau\zeta - \zeta\rho|$. If $\rho = \infty,\tau,\zeta \in \nR$, it simplifies to $|\zeta-\tau|$. If $\tau = \infty,\rho,\zeta \in \nR$, it simplifies to $|\rho+\zeta -1|$. This simplified forms of the ratio of the similarity follow from the identities (\ref{eqthreeb}), (\ref{eqtwoinfty}) and (\ref{eqtwoone}).
\end{proof}

\begin{corollary} \label{ctds}
Let $\rho, \tau, \zeta \in \nRc$ and let $T$ be an oriented non-equilateral triangle. Then $\CT_\tau \bigl(\CT_\rho(T)\bigr)$ is directly similar to $\CT_\zeta(T)$ if and only if
$\zeta\Bo\rho \in  \tau \Bo\varkappa_{_T} \Bo {\mathbb S}$.
\end{corollary}
\begin{proof}
By Proposition~\ref{prepeat} and Corollary~\ref{rho-T-V}, the triangle $\CT_\tau \bigl(\CT_\rho(T)\bigr)$ is directly similar to $\CT_\zeta(T)$ if and only if $\CT_\zeta\bigl(\CT_\tau \bigl(\CT_\rho(T)\bigr)\bigr)$ is directly similar to $T$. By Theorem~\ref{reduction 01} and Lemma~\ref{lcomp}, the last statement is equivalent to $\rho \Bo\Ton{\tau}\Bo\zeta \in \varkappa_{_T}\Bo\nS$. The statement follows by applying the operation $\Bo\tau$ on the both sides of the last relation.
\end{proof}

\begin{theorem} \label{tinsim}
Let  $\rho, \tau, \zeta \in \nRc$ and let $T$ be an oriented non-equilateral triangle. Then $\CT_\tau \bigl(\CT_\rho(T)\bigr)$ is reversely similar to $\CT_\zeta(T)$ if and only if $\tau \Bo \zeta \in \rho\Bo\nS$. In particular, $\CT_{\tau}(T)$ and $\CT_{\zeta}(T)$ are reversely similar if and only if $\tau \Bo \zeta \in\varkappa_{_T} \Bo\nS$.
\end{theorem}
\begin{proof}
The triangle $\CT_\tau \bigl(\CT_\rho(T)\bigr)$ is reversely similar to $\CT_\zeta(T)$ if and only if there exists $\sigma \in \nS$ and $\lambda\in \nR$ such that
$\Vs{M}_\sigma \Vs{M}_\tau \Vs{M}_\rho \Vb{t} = \lambda \Vs{M}_\zeta \Vb{t}$. By Theorem~\ref{tcomp} and Proposition~\ref{prsu}, the last equality is equivalent to $\lambda = 1$ and  $\rho\Bo\Ton{\tau}\Bo\sigma = \zeta$. Recalling that $\sigma\in\nS$, the last statement is equivalent to $\tau \Bo\zeta \in \rho\Bo\nS$. Setting $\rho = \varkappa_{_T}$ yields the special case.
\end{proof}

\section{Isosceles triangles}

Let $T$ be an oriented non-equilateral triangle. In this section we identify the isosceles triangles in the family $\CT_\rho(T), \rho \in \nRc$. As we will see, these triangles play an important role in this family. Recall that $\nT = \{0, 1/2, 1, 2, \infty, -1\}$.

\begin{proposition} \label{pit}
Let $T$ be an oriented  non-equilateral triangle and $\rho \in \nRc$. The triangle $\CT_\rho(T)$ is an isosceles triangle if and only if $\rho \in  \raisebox{0pt}{$\sqrt[\Bo]{\varkappa_{_T}}$}\Bo \nT$.
Let $\rho \in \sqrt[\Bo]{\varkappa_{_T}}\Bo\nS$ and $\zeta \in  \sqrt[\Bo]{\varkappa_{_T}}\Bo(\nT\setminus \nS)$. Then $\CT_\rho(T)$ is wide (narrow) if and only if $\CT_\zeta(T)$ is narrow (wide).
\end{proposition}
\begin{proof}
It follows from Remark~\ref{risop} that the triangle $\CT_\rho(T)$ is isosceles if and only if there exists $\sigma\in\nS$ such that $\Vs{M}_\sigma \Vs{M}_\rho \Vb{t} = \Vs{M}_\rho \Vb{t}$, or, equivalently $\Vs{M}_\rho \Vs{M}_\sigma \Vs{M}_\rho \Vb{t} = \Vb{t}$. By Theorems~\ref{tcomp} and~\ref{pgamma}, the last equality holds if and only if $\rho\Bo\rho\Bo\Ton{\sigma} = \varkappa_{_T}$. If $\sigma = 0$, the solutions of the last equation are $\sqrt[\Bo]{\varkappa_{_T}} \Bo\{0,2\}$.
If $\sigma = \infty$, the solutions are $\sqrt[\Bo]{\varkappa_{_T}} \Bo\{1,-1\}$.
If $\sigma = 1$, the solutions are $\sqrt[\Bo]{\varkappa_{_T}} \Bo\{\infty,\tfrac{1}{2}\}$. This proves the first claim.

To prove the second claim let $\rho \in \sqrt[\Bo]{\varkappa_{_T}}\Bo\nS$ and $\zeta \in \sqrt[\Bo]{\varkappa_{_T}}\Bo(\nT\setminus \nS)$. By Theorem~\ref{reduction 01}, $\CT_\rho(T)$ is similar to $\CT_{\rho_1}(T)$ for any $\rho_1 \in \sqrt[\Bo]{\varkappa_{_T}}\Bo\nS$. Therefore, without loss of generality, we can assume that $\zeta = \rho \Bo 2$. Thus, $\zeta =  0 \Bo 2 \Bo \rho$ and
\begin{equation*}
\Vs{M}_\zeta = \Vs{M}_0 \Vs{M}_2 \Vs{M}_\rho =
\Vs{Q}\Rot\bigl(2(\Phi^{-1}(0) - \Phi^{-1}(2))\bigr)\Vs{Q}^\top
\Vs{M}_\rho.
\end{equation*}
Since $2(\Phi^{-1}(0) - \Phi^{-1}(2)) = -\pi$, if $\Vs{M}_\rho \Vb{t}$ is wide, then $\Vs{M}_\zeta \Vb{t}$ is narrow, and conversely.
The proposition is proved.
\end{proof}

Since $\nT\setminus\nS = \tfrac{1}{2} \Bo\, \nS$, exactly one element in each of the sets $\sqrt[\Bo]{\varkappa_{_T}}\Bo\nS$ and
\mbox{$\sqrt[\Bo]{\varkappa_{_T}}\Bo(\nT\setminus \nS)$} belongs to $\nI$. Proposition~\ref{pimath} implies that those special elements are $\sqrt[\Bo]{\pI(\varkappa_{_T}\!)}$ and \mbox{$\tfrac{1}{2} \Bo \sqrt[\Bo]{\pI(\varkappa_{_T}\!)}$.} These observations together with Proposition~\ref{pit} prove the first statement of the next theorem.

\begin{figure}

\setlength{\abovecaptionskip}{4pt}%
\setlength{\belowcaptionskip}{-8pt}%

\noindent\hspace*{-5pt}\begin{minipage}{0.5\textwidth}

\psfrag{mu}[][]{\begin{picture}(0,0)
            \put(-4,-14){\makebox(0,0)[l]{$\mu_{_T}$}}
                        \end{picture}}

\psfrag{nu}[][]{\begin{picture}(0,0)
            \put(-4,-14){\makebox(0,0)[l]{$\nu_{_T}$}}
                        \end{picture}}

\psfrag{zr}[][]{\begin{picture}(0,0)
            \put(0,-1){\makebox(0,0)[l]{$z_\rho^2$}}
                        \end{picture}}

\psfrag{yr}[][]{\begin{picture}(0,0)
            \put(0,-10){\makebox(0,0)[l]{$y_\rho^2$}}
                        \end{picture}}

\psfrag{xr}[][]{\begin{picture}(0,0)
            \put(0,0){\makebox(0,0)[l]{$x_\rho^2$}}
                        \end{picture}}

\psfrag{z}[][]{\begin{picture}(0,0)
            \put(-3,2){\makebox(0,0)[l]{$0$}}
                        \end{picture}}
\psfrag{h}[][]{\begin{picture}(0,0)
            \put(-3,0){\makebox(0,0)[l]{$\tfrac{1}{2}$}}
                        \end{picture}}
\psfrag{o}[][]{\begin{picture}(0,0)
            \put(-3,2){\makebox(0,0)[l]{$1$}}
                        \end{picture}}

\resizebox{!}{!}{%
  \includegraphics{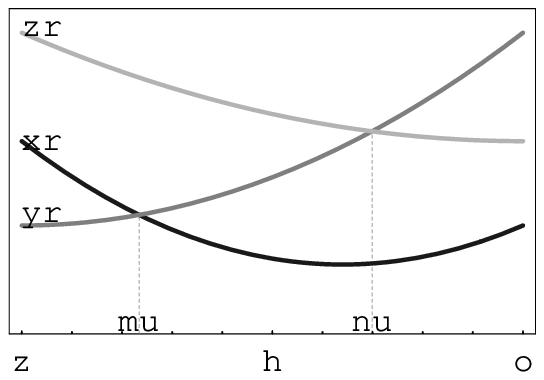}}
    \caption{An increasing $T$}
\label{fincT}
\end{minipage}\hspace*{-5pt}
\begin{minipage}{0.5\textwidth}

\psfrag{mu}[][]{\begin{picture}(0,0)
            \put(-4,-14){\makebox(0,0)[l]{$\mu_{_T}$}}
                        \end{picture}}

\psfrag{nu}[][]{\begin{picture}(0,0)
            \put(-4,-14){\makebox(0,0)[l]{$\nu_{_T}$}}
                        \end{picture}}

\psfrag{zr}[][]{\begin{picture}(0,0)
            \put(2,0){\makebox(0,0)[l]{$z_\rho^2$}}
                        \end{picture}}

\psfrag{yr}[][]{\begin{picture}(0,0)
            \put(2,-1){\makebox(0,0)[l]{$y_\rho^2$}}
                        \end{picture}}

\psfrag{xr}[][]{\begin{picture}(0,0)
            \put(2,4){\makebox(0,0)[l]{$x_\rho^2$}}
                        \end{picture}}

\psfrag{z}[][]{\begin{picture}(0,0)
            \put(-3,2){\makebox(0,0)[l]{$0$}}
                        \end{picture}}
\psfrag{h}[][]{\begin{picture}(0,0)
            \put(-3,0){\makebox(0,0)[l]{$\tfrac{1}{2}$}}
                        \end{picture}}
\psfrag{o}[][]{\begin{picture}(0,0)
            \put(-3,2){\makebox(0,0)[l]{$1$}}
                        \end{picture}}

\resizebox{!}{!}{%
  \includegraphics{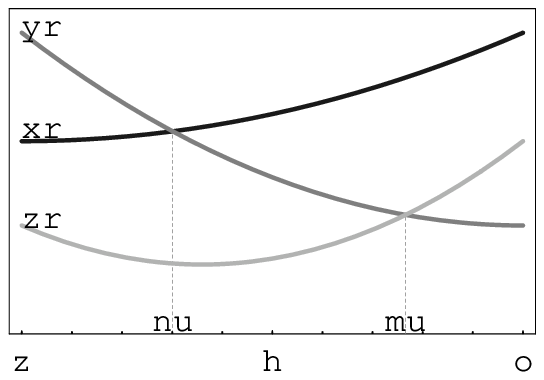}}
    \caption{A decreasing $T$}
\label{fdecT}
\end{minipage}

\end{figure}

\begin{theorem} \label{tit}
Let $T$ be an oriented non-equilateral triangle.
\begin{enumerate}[{\rm (a)}]
\item \label{titia}
There exist unique numbers  $\mu_{_T}, \nu_{_T} \in \nI$ such that $\CT_{\mu_{_T}}(T)$ is wide and $\CT_{\nu_{_T}}(T)$ is narrow.
\item \label{titib}
If $T$ is wide, $\mu_{_T} = \pI(\varkappa_{_T}\!) = 0$ and $\nu_{_T}
= 1/2$. If $T$ is narrow, $\nu_{_T} = \pI(\varkappa_{_T}\!) =  0$
and $\mu_{_T} = 1/2$.
\item \label{titic}
If $T$ is increasing then $\mu_{_T} = \sqrt[\Bo]{\pI(\varkappa_{_T}\!)}, \nu_{_T} = \tfrac{1}{2} \BoI \sqrt[\Bo]{\pI(\varkappa_{_T}\!)}$, and $\mu_{_T} < \varkappa_{_T} < \nu_{_T}$.
\item \label{titid}
If $T$ is decreasing then
$\nu_{_T} = \sqrt[\Bo]{\pI(\varkappa_{_T}\!)}, \mu_{_T}=\tfrac{1}{2} \BoI \sqrt[\Bo]{\pI(\varkappa_{_T}\!)}$, and $\nu_{_T} < \varkappa_{_T} < \mu_{_T}$.
\item \label{titie}
A triangle $\CT_{\rho}(T)$ is wide if and only if $\rho \in \mu_{_T} \Bo \nS$.
\item \label{titif}
A triangle $\CT_{\rho}(T)$ is narrow if and only if $\rho \in \nu_{_T} \Bo \nS$.
\end{enumerate}
\end{theorem}
\begin{proof}
If $T$ is wide, then clearly $\CT_0(T)$ is wide and $\CT_{1/2}(T)$ is narrow.  Similarly, if $T$ is narrow, then $\CT_0(T)$ is narrow and $\CT_{1/2}(T)$ is wide. Hence, (\ref{titib}) holds.

To prove (\ref{titic}), let $T$ be an increasing oriented triangle. Without loss of generality, assume that $a < b < c$ are its sides.  With the notation introduced in \eqref{Crho}, \eqref{Brho}, \eqref{Arho}, we immediately see that, for all $\rho \in \mathbb I$, $x_\rho^2 < z_\rho^2$ (see Figure~\ref{fincT}). Also,
\[
y_0^2 = a^2 < x_0^2 = b^2 \qquad \text{and} \qquad y_1^2 = c^2 > x_1^2 = a^2.
\]
Therefore, since $x_\rho$ and $y_\rho$ are continuous functions of $\rho$, there exists $\rho_1 \in \nI$ such that $x_{\rho_1} = y_{\rho_1}$. Since $z_{\rho_1} > y_{\rho_1} = x_{\rho_1}$, $\CT_{\rho_1}(T)$ is wide.
Similarly,
\[
y_0^2 = a^2 < z_0 = c^2 \qquad \text{and} \qquad y_1^2 = c^2 > z_1^2 = b^2.
\]
Therefore, there exists $\rho_2 \in \nI$ such that $z_{\rho_2} = y_{\rho_2}$. Since $x_{\rho_2} < y_{\rho_2} = z_{\rho_2}$, $\CT_{\rho_2}(T)$ is narrow.
Since $x_\rho^2 < z_\rho^2$ for all $\rho \in \nI$ we have $\rho_1 < \rho_2$. As $\sqrt[\Bo]{\pI(\varkappa_{_T}\!)}$ and $\tfrac{1}{2} \BoI \sqrt[\Bo]{\pI(\varkappa_{_T}\!)}$ are the only values of $\rho \in \nI$ for which $\CT_\rho(T)$ is isosceles and since we have  $\sqrt[\Bo]{\pI(\varkappa_{_T}\!)} < \tfrac{1}{2} \BoI \sqrt[\Bo]{\pI(\varkappa_{_T}\!)}$, it follows that
\[
\rho_1 = \sqrt[\Bo]{\pI(\varkappa_{_T}\!)} \qquad \text{and} \qquad \rho_2 = \tfrac{1}{2} \BoI \sqrt[\Bo]{\pI(\varkappa_{_T}\!)}.
\]
This proves (\ref{titic}). Analogous reasoning proves (\ref{titid}), see Figure~\ref{fdecT}. The items (\ref{titie}) and (\ref{titif}) follow from Theorem~\ref{reduction 01}.
\end{proof}

\begin{remark} \label{rmuf}
As it was pointed out in Section~\ref{sbn}, the sides of an arbitrary  oriented  non-equilateral triangle $T$ can be labeled counterclockwise uniquely in such a way that $a \leq b < c$ or $a \geq b > c$. Then, $\varkappa_{_T} = (a^2-b^2)/(a^2-c^2)  \in \nI$ and $\mu_{_T} = \sqrt[\Bo]{\varkappa_{_T}}$ if $T$ is increasing and $\mu_{_T} = \tfrac{1}{2}\BoI\sqrt[\Bo]{\varkappa_{_T}}$ if $T$ is decreasing. Simplifying the corresponding formulas, $\mu_{_T}$ can now be expressed in terms of $a,b,c$ as
\[
\mu_{_T} = \left\{ \begin{array}{cc}
\dfrac{b^2-a^2}{c^2-a^2+\sqrt{a^4+b^4+c^4-a^2b^2-b^2c^2-c^2a^2}} & \text{if} \quad a \leq b < c, \\[14pt]
\dfrac{a^2-c^2}{b^2-c^2+\sqrt{a^4+b^4+c^4-a^2b^2-b^2c^2-c^2a^2}} & \text{if} \quad a \geq b > c.
\end{array} \right.
\]
Notice that
\[
a^4+b^4+c^4-a^2b^2-b^2c^2-c^2a^2 = \frac{1}{2} \bigl((a^2-b^2)^2 + (b^2-c^2)^2 +(c^2-a^2)^2\bigr).
\]
This shows that the quantity under the square root is nonnegative and that the value of $\mu_{_T}$ for a narrow triangle $a = b > c$ is $1/2$; see also Theorem \ref{tit} (b).

As shown in Remark~\ref{rappI}, the groups $(\nI,\BoI)$ and $(\nI,\oplus)$ are very close to each other. Therefore, a very good approximation for $\sqrt[\Bo]{\varkappa_{_T}}$ is $\varkappa_{_T}/2$ and a very good approximation for $\tfrac{1}{2}\BoI\sqrt[\Bo]{\varkappa_{_T}}$ is $(1+\varkappa_{_T})/2$. In fact, with these approximations, the absolute error is, in both cases, less than $0.021$, while the relative error in the first case is less than $0.067$ and  less than $0.033$ in the second case. Hence, the above long formulas for $\mu_{_T}$ can be well approximated by
\[
\frac{b^2-a^2}{2(c^2-a^2)} \quad \text{if} \quad a \leq b < c \quad \text{and} \quad \frac{2a^2 - b^2 - c^2}{2(a^2-c^2)} \quad \text{if} \quad a \geq b > c.
\]
\end{remark}


Define
\[
{\mathbb M}_{T} := \bigl[\min\{\mu_{_T}, \nu_{_T}\},\max\{\mu_{_T}, \nu_{_T}\} \bigr] = \bigl[\sqrt[\Bo]{\pI(\varkappa_{_T}\!)}, \tfrac{1}{2}\BoI\sqrt[\Bo]{\pI(\varkappa_{_T}\!)}\bigr].
\]

Recall that in Theorem~\ref{reduction 01} we showed that all triangles in the family  $\CT_{\mathbb R}(T) \cup \{\CT_\infty(T)\}$, up to direct similarity, can be found in $\CT_{{\mathbb I}} (T)$. In the next theorem, we show that all triangles from $\CT_{\mathbb I}(T)$, up to similarity, can be found in $\CT_{{\mathbb M}_{T}} (T)$.

\begin{theorem}[Reduction to $\CT_{{\mathbb M}_{T}} (T)$] \label{CJT}
Let $T$ be an oriented non-equilateral triangle. The mapping $\imath: \nI \to \nI$ defined by
\begin{equation} \label{involution}
\imath(\xi)= {\mathsf p}(\varkappa_{_T})\!\BoI \TonI{\xi}, \quad \xi \in \nI,
\end{equation}
maps the interior of $\mathbb M_{\;T}$ onto its exterior ${\mathbb I} \setminus\! {\mathbb M}_{\;T}$. The triangles $\CT_{\xi} (T)$ and $\CT_{\imath(\xi)} (T)$ are reversely similar for all $\xi \in {\mathbb M}_{\;T}$.
\end{theorem}

\begin{proof}
The first claim follows from Proposition~\ref{pimath}~\!(\ref{pimathid}). To prove the inverse similarity of the triangles $\CT_\xi (T)$ and $\CT_{\imath(\xi)} (T)$, we use the last claim of Theorem~\ref{tinsim} with $\tau = \imath(\xi)$ and $\zeta = \xi$. Then, the condition $\tau\Bo \zeta = (\imath(\xi)) \Bo \xi \in \varkappa_{_T} \Bo \nS$ is equivalent to $\pI\bigl((\imath(\xi)) \Bo \xi \bigr) = \pI(\varkappa_{_T}\!)$, which, in turn, is equivalent to $(\imath(\xi)) \BoI \xi = \pI(\varkappa_{_T}\!)$. Since the last equality is trivial, the theorem is proved.
\end{proof}

\begin{remark} \label{rsameo}
Let $T$ be scalene triangle. Then $\pI(\varkappa_{_T}\!) \neq 0$ and thus, by \eqref{eqfps}, $1/2$ is an interior point of ${\mathbb M}_{T}$. It follows from the proof of Theorem~\ref{tit}, see also Figures~\ref{fincT} and~\ref{fdecT}, that if $T$ is increasing, then all the triangles in the interior of ${\mathbb M}_{T}$ are increasing, and for a decreasing $T$ all the triangles in the interior of ${\mathbb M}_{T}$ are decreasing. Therefore the median triangle of an arbitrary scalene triangle has the same orientation as the host triangle.
\end{remark}

\section{Related shape functions} \label{sshape}

Let $T=(a,b,c)$ be a given oriented triangle. The Brocard angle $\omega_{_T}$ and, hence, the cone angle $\gamma_{_T}$ are related to the following shape
function (that is, a complex-valued function that characterizes the
similarity of two geometric objects) introduced by Hajja in
\cite{Hajja2}:
\[
\sigma (T)=\frac{a^2+e^{-2\pi i/3}b^2+e^{2\pi
i/3}c^2}{a^2+b^2+c^2}.
\]
Using \cite[Theorem~2.4]{Hajja2} and \eqref{eqbaca} we get $
|\sigma(T)| =  (\tan \gamma_{_T})/\sqrt{2}.$
As already observed in
\cite[Theorem~3.1~(e)]{Hajja2} or
\cite[Theorem~5.1~(2c)]{Hajja1}, the equality of the Brocard angles of two triangles $T$ and $V$ implies the existence of $\rho \in \nR$ such
that $T$ and $\CT_\rho (V)$ (${\mathcal H}_\rho (V)$ in \cite{Hajja1,
Hajja2})  are similar. In the following theorem we
extend this statement with several equivalences
that relate to the special values $\mu_{_T}$ and $\nu_{_T}$. These equivalences allow us, for example, to identify exactly the parameter $\rho$ for which $T$ is similar to $\CT_\rho (V)$. This leads to several shape functions such as, $T\mapsto \omega_{_T} + i\mu_{_T}$, $T\mapsto \gamma_{_T} + i\mu_{_T}$, $T\mapsto \omega_{_T} + i\nu_{_T}$  and  $T \mapsto \gamma_{_T} + i\nu_{_T}$.

\begin{theorem} \label{th1}
Let $T$ and $V$ be oriented triangles. The following statements are equivalent.
\begin{enumerate}[{\rm (a)}]
\item \label{thia0}
There exist $\xi, \zeta \in \nRc$ such that $\CT_\zeta(V)$ and $\CT_\xi(T)$ are similar.
\item \label{thib}
There exists $\rho \in {\mathbb M}{_T}$ such that $V$ and $\CT_\rho(T)$ are similar.
\item \label{thia}
There exists $\tau \in {\mathbb M}{_V}$ such that $T$ and $\CT_\tau(V)$ are similar.
\item \label{thic}
$\CT_{\mu_{_{T}}}(T)$ is similar to $\CT_{\mu_{_{V}}}(V)$.
\item \label{thid}
$\CT_{\nu_{_{T}}}(T)$ is similar to $\CT_{\nu_{_{V}}}(V)$.
\item \label{thie}
$\gamma_{_T} = \gamma_{_V}$.
\item \label{thif}
$\omega_{_T} = \omega_{_V}$.
\end{enumerate}
\end{theorem}
\begin{proof}
We will prove (\ref{thia0})$\Rightarrow$(\ref{thib})$\Rightarrow$(\ref{thia})
$\Rightarrow$(\ref{thic})$\Rightarrow$(\ref{thia0}), (\ref{thib})$\Leftrightarrow$(\ref{thie}) and (\ref{thic})$\Leftrightarrow$(\ref{thid}); (\ref{thie})$\Leftrightarrow$(\ref{thif}) follows from Proposition~\ref{pbaca}.

Assume (\ref{thia0}). Then Corollary~\ref{rho-T-V} yields that  $\CT_\zeta\bigl(\CT_\zeta(V)\bigr)$ and  $\CT_\zeta\bigl(\CT_\xi(T)\bigr)$ are similar.  By Proposition~\ref{prepeat},  $\CT_\zeta\bigl(\CT_\zeta(V)\bigr)$ is directly similar to $V$ and, by Theorem~\ref{tinsim},
$\CT_\zeta\bigl(\CT_\xi(T)\bigr)$ is reversely similar to $\CT_{\xi\Bo\Ton{\zeta}}(T)$. Hence, $V$ is similar to $\CT_{\xi\Bo\Ton{\zeta}}(T)$. By Theorems~\ref{reduction 01}
and~\ref{CJT}, (\ref{thib}) follows.  Assume (\ref{thib}). By applying $\CT_\rho$ we get that $\CT_\rho(V)$ is similar to $\CT_\rho(\CT_\rho(T))$. By Proposition~\ref{prepeat} and Theorems~\ref{reduction 01} and~\ref{CJT}, (\ref{thia}) follows.
Now assume (\ref{thia}). Then $\CT_{\mu_{_{T}}}(T)$ is similar to $\CT_{\mu_{_{T}}}\bigl(\CT_\tau(V)\bigr)$. By Theorem~\ref{tinsim},
$\CT_{\mu_{_{T}}}(T)$ is similar to $\CT_{\tau\Bo\Ton{\mu_{_{T}}}}(V)$.
Theorem~\ref{tit}~(\ref{titie}) yields (\ref{thic}).
Since (\ref{thic}) is a special case of (\ref{thia0}), the first sequence of implications is proved.

The statement (\ref{thib}) is equivalent to $\Vs{M}_\rho \Vb{t} = \lambda \Vb{v}$. By Theorem~\ref{pgamma}, the last equality is equivalent to $(\Vb{t}\cdot\Vb{r})/\|\Vb{t}\| = (\Vb{v}\cdot\Vb{r})/\|\Vb{v}\|$. Since
$(\Vb{t}\cdot\Vb{r})/\|\Vb{t}\| = - \cos \gamma_{_T}$, the equivalence of (\ref{thib}) and (\ref{thie}) is proved.

Finally, we prove (\ref{thic})$\Leftrightarrow$(\ref{thid}). Since $\sqrt[\Bo]{\pI(\varkappa_{_T}\!)}\Bo \tfrac{1}{2}\Bo \sqrt[\Bo]{\pI(\varkappa_{_T}\!)} \Bo \Ton{\bigl(\tfrac{1}{2}\bigr)} = \pI(\varkappa_{_T}\!)$, we have $\mu_{_T} \Bo \nu_{_T} \Bo \Ton{\bigl(\tfrac{1}{2}\bigr)} = \pI(\varkappa_{_T}\!)$, or equivalently $\mu_{_T}\Bo \nu_{_T} \Bo \Ton{\bigl(\tfrac{1}{2}\bigr)} \in \varkappa_{_T} \Bo \nS$. Theorems~\ref{reduction 01} and~\ref{CJT} together with  Lemma~\ref{lcomp} now yield that  $\CT_\nu\bigl(\CT_{1/2}\bigl(\CT_\mu(T)\bigr)\bigr)$ similar to $T$. Consequently, $\CT_{1/2}\bigl(\CT_{\mu_{_T}}(T)\bigr)$ is similar to $\CT_{\nu_{_T}}(T)$. With the analogous similarity for $V$, we have that the statement (\ref{thid}) is equivalent to $\CT_{1/2}\bigl(\CT_{\mu_{_T}}(T)\bigr)$ being similar to $\CT_{1/2}\bigl(\CT_{\mu_{_V}}(V)\bigr)$, and this similarity is equivalent to the statement (\ref{thic}).
\end{proof}

We have pointed out in Remark~\ref{rsameo} that, for an increasing $T$, all the triangles $\CT_\rho(T)$ with $\rho$ in the interior of $\nM_T$ are increasing. The reasoning from the proof of Theorem~\ref{tit} also yields that, for a wide triangle $W$, all the triangles $\CT_\rho(W)$ with $\rho$ in the interior of $\nM_W$ are increasing. Similarly, for a narrow $N$, all the triangles $\CT_\rho(N)$ with $\rho$ in the interior of $\nM_N$ are decreasing. Therefore, in the next proposition, the wide triangles are included in the family of increasing triangles, while the narrow triangles are included in the family of decreasing triangles.

\begin{proposition} \label{cbijection}
Assume that any of the equivalent conditions in Theorem~{\rm~\ref{th1}} is  satisfied. For $\xi \in \nM_T$, set
 \[
Z(\xi) = \begin{cases}
 \mu_{_V}\!\!\BoI\TonI{\mu_{_T}}\!\BoI\xi &
   \text{if \ $T$ and $V$ are both increasing or both decreasing}, \\[8pt]
 \mu_{_V}\!\!\BoI \mu_{_T}\!\BoI\TonI{\xi}  &
  \text{\begin{minipage}[t]{2.7in}
 if \ $T$ is increasing and $V$ is decreasing or \\
 \rule{10pt}{0pt} $T$ is decreasing and $V$ is increasing.
 \end{minipage}}
\end{cases}
 \]
Then $Z:\nM_T \to \nM_V$ is a monotonic bijection. If $T$ and $V$ are both increasing or both decreasing, then $\CT_{Z(\xi)}(V)$ and $\CT_\xi(T)$ are directly similar for all $\xi \in \nM_T$. If $T$ is increasing and $V$ is decreasing or if $T$ is decreasing and $V$ is increasing, then $\CT_{Z(\xi)}(V)$ and $\CT_\xi(T)$ are reversely similar for all $\xi \in \nM_T$.
\end{proposition}
\begin{proof}
If $T$ and $V$ are both increasing, then $\nu_{_T} = \tfrac{1}{2}\BoI\mu_{_T} > \mu_{_T}$ and $\nu_{_V} = \tfrac{1}{2}\BoI\mu_{_V} > \mu_{_V}$. If $T$ and $V$ are both decreasing,  then $\mu_{_T} = \tfrac{1}{2}\BoI\nu_{_T} > \nu_{_T}$ and $\mu_{_V} = \tfrac{1}{2}\BoI\nu_{_V} > \nu_{_V}$. In either case
\( Z(\mu_{_T}) = \mu_{_V}\)  and \(Z(\nu_{_T}) = \nu_{_V}.\) By Proposition~\ref{padx}, $Z$ maps $\nM_T$ onto $\nM_V$ as an increasing bijection. Since $\CT_{\mu_{_T}}(T)$ is directly similar to $\CT_{\mu_{_V}}(V)$, by Proposition~\ref{prepeat} $T$ is directly similar to $\CT_{\mu_{_T}}\bigl(\CT_{\mu_{_V}}(V)\bigr)$. Therefore, $\CT_{\xi}(T)$ is directly similar to $\CT_{\xi}\bigl(\CT_{\mu_{_T}}\bigl(\CT_{\mu_{_V}}(V)\bigr)\bigr)$, and hence,  by Theorem~\ref{tcomp},
$\CT_\xi(T)$ is directly similar to $\CT_{Z(\xi)}(V)$.

If $T$ is decreasing and $V$ is increasing, then $\mu_{_T} = \tfrac{1}{2}\BoI\nu_{_T} > \nu_{_T}$ and $\nu_{_V} = \tfrac{1}{2}\BoI\mu_{_V} > \mu_{_V}$. Then, again, \( Z(\mu_{_T}) = \mu_{_V}\) and \(Z(\nu_{_T}) = \nu_{_V}.\)  These equalities hold if $T$ is increasing and $V$ is decreasing, as well. Proposition~\ref{pimath} implies that $Z$ maps $\nM_T$ onto $\nM_V$ as a decreasing bijection. Now, we prove the reverse similarity of triangles $\CT_{Z(\xi)}(V)$ and $\CT_\xi(T)$.  Let $T^\prime$ be a triangle which is reversely congruent to $T$. Then $T^\prime$ is increasing and, by Proposition~\ref{pinvcon} and Theorem~\ref{tit}~\!(\ref{titia}),
$\mu_{_{T^\prime}} = 1- \mu_{_T} = \TonI{\mu_{_T}}$. In general, $\CT_{1-\xi}(T^\prime)$ is reversely congruent to $\CT_\xi(T)$. Since $T^\prime$ is increasing, the first part of this proof yields that $\CT_{1-\xi}(T^\prime)$ is directly similar to $\CT_\zeta(V)$ for $\zeta = \mu_{_V}\BoI \TonI{\mu_{_{T^\prime}}} \BoI (1-\xi)$. Since $\mu_{_{T^\prime}} = \TonI{\mu_{_T}}$ and $1-\xi = \TonI{\xi}$, $\CT_{1-\xi}(T^\prime)$ is directly similar to $\CT_{Z(\xi)}(V)$. As $\CT_{1-\xi}(T^\prime)$ is reversely congruent to $\CT_\xi(T)$, we get that $\CT_\xi(T)$ is reversely similar to $\CT_{Z(\xi)}(V)$.
\end{proof}

\begin{remark} \label{rVsT}
If any of the conditions of Theorem~{\rm~\ref{th1}} hold, we can use Proposition~\ref{cbijection} to obtain a formula for $\tau \in \nM_{V}$ such that $T$ is similar to $\CT_\tau(V)$. If $T$ and $V$ have the same orientation, then
\[
\tau = Z(\pI(\varkappa_{_T}\!)) = \mu_{_V}\!\BoI\TonI{\mu_{_T}}\!\BoI \pI(\varkappa_{_T}\!) = \mu_{_V}\!\BoI\mu_{_T},
\]
and $T$ is directly similar to $\CT_\tau(V)$. If $T$ and $V$ have opposite orientation, then
\[
\tau = Z(\pI(\varkappa_{_T}\!)) = \mu_{_V}\!\BoI\mu_{_T}\!\BoI \TonI{(\pI(\varkappa_{_T}\!))} = \mu_{_V}\!\BoI\TonI{\mu_{_T}},
\]
and $T$ is reversely similar to $\CT_\tau(V)$. Moreover, in this case we have that $T$ is directly similar to $\CT_0(\CT_\tau(V))$. By Corollary~\ref{ctds}, $\CT_0(\CT_\tau(V))$ is directly similar to $\CT_\zeta(V)$ with $\zeta \Bo \mu_{_V}\BoI\TonI{\mu_{_T}} \in \varkappa_{_T} \Bo \nS$. To find such $\zeta$, we apply $\pI$ to both sides of the last membership to get the equation $\zeta\BoI\mu_{_V}\!\BoI\TonI{\mu_{_T}} = \pI(\varkappa_{_T}\!)$, whose solution is  $\zeta =  \TonI{\mu_{_V}}\BoI\mu_{_T}\!\BoI \pI(\varkappa_{_T}\!) = \mu_{_V}\!\BoI\mu_{_T}$.

Thus, amazingly, for arbitrary non-equilateral oriented triangles $T$ and $V$ with the same Brocard angle $T$ is directly similar to $\CT_{\mu_{_V}\!\!\BoI\mu_{_T}}\!(V)$. It is worth noting that when $T$ and $V$ have opposite orientation we have $\mu_{_V}\!\!\BoI\mu_{_T}\! \notin \nM_V$. However, by definition of $\BoI$, always $\mu_{_V}\!\BoI\mu_{_T}\! \in \nI$.
\end{remark}

\begin{example} \label{exVsT}
Here we illustrate the last claim in Remark~\ref{rVsT} using the triangles \begin{align*}
T & = \left(
  \sqrt{5\sqrt{7}+5}, \sqrt{5\sqrt{7} - 4}, \sqrt{5 \sqrt{7} -1}\right), \\
V & =  \left(
      2 \sqrt{5 \sqrt{7}+4},2\sqrt{5\sqrt{7}+1},2 \sqrt{5\sqrt{7} - 5}
      \right).
\end{align*}
The triangle $T$ is increasing and $V$ is decreasing. We calculate $\tan \omega_{_T} = \tan \omega_{_V} = \sqrt{5}/7$,  $\mu_{_T} = \bigl(3-\sqrt{7}\bigr)/2$, $\mu_{_V} = \sqrt{7}-5$ and $\mu_{_V}\!\BoI\mu_{_T} = 4/5$. In Figure~\ref{fSBAa} the reader can see that $\CT_{4/5}(T) = XY\!Z$ is directly similar to $V=KLM$, while in Figure~\ref{fSBAb} the reader can see that $\CT_{4/5}(V) = XY\!Z$ is directly similar to $T=ABC$. Since $T$ is increasing and $V$ is decreasing, Remark~\ref{rsameo} yields that $4/5 \not\in \nM_{T}$ and $4/5 \not\in \nM_{V}$. This is easy to verify since in this case $\mu_{_T} = \nu_{_V} = \bigl(3-\sqrt{7}\bigr)/2 \approx 0.1771$, $\nu_{_T} = \mu_{_V} = \sqrt{7}-5 \approx 0.6458$.

\begin{figure}[ht]

\psfrag{X}[][]{\begin{picture}(0,0)
            \put(-4,0){\makebox(0,0)[l]{$X$}}
                        \end{picture}}
\psfrag{Y}[][]{\begin{picture}(0,0)
            \put(-3,0){\makebox(0,0)[l]{$Y$}}
                        \end{picture}}
\psfrag{Z}[][]{\begin{picture}(0,0)
            \put(-7,-7){\makebox(0,0)[l]{$Z$}}
                        \end{picture}}

\psfrag{A}[][]{\begin{picture}(0,0)
            \put(-10,2){\makebox(0,0)[l]{$A$}}
                        \end{picture}}
\psfrag{B}[][]{\begin{picture}(0,0)
            \put(-2,-3){\makebox(0,0)[l]{$B$}}
                        \end{picture}}
\psfrag{C}[][]{\begin{picture}(0,0)
            \put(-6,-3){\makebox(0,0)[l]{$C$}}
                        \end{picture}}
\psfrag{App}[][]{\begin{picture}(0,0)
            \put(11,2){\makebox(0,0)[l]{$A^\prime$}}
                        \end{picture}}

\psfrag{Ap}[][]{\begin{picture}(0,0)
            \put(-8,-3){\makebox(0,0)[l]{$A_\rho$}}
                        \end{picture}}

\psfrag{Bp}[][]{\begin{picture}(0,0)
            \put(-4,-13){\makebox(0,0)[l]{$B_\rho$}}
                        \end{picture}}

\psfrag{Cp}[][]{\begin{picture}(0,0)
            \put(-23,-2){\makebox(0,0)[l]{$C_\rho$}}
                        \end{picture}}

\psfrag{K}[][]{\begin{picture}(0,0)
            \put(-2,0){\makebox(0,0)[l]{$K$}}
                        \end{picture}}

\psfrag{L}[][]{\begin{picture}(0,0)
            \put(-6,-2){\makebox(0,0)[l]{$L$}}
                        \end{picture}}

\psfrag{M}[][]{\begin{picture}(0,0)
            \put(-2,-4){\makebox(0,0)[l]{$M$}}
                        \end{picture}}

\resizebox{!}{!}{%
  \includegraphics{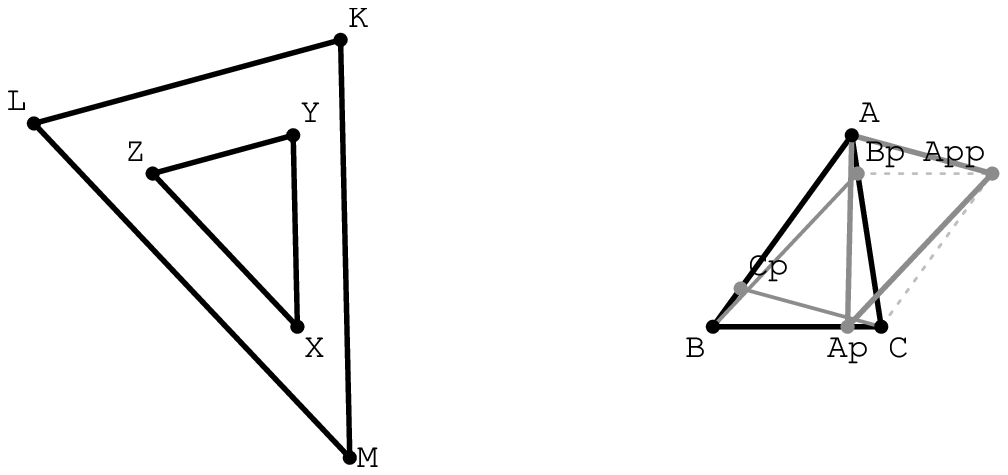}}
    \caption{$T=ABC, V=KLM, \rho =4/5, \CT_{\rho} (T) = XY\!Z$ is directly similar to $V$}
\label{fSBAa}

\psfrag{Kp}[][]{\begin{picture}(0,0)
            \put(-12,-1){\makebox(0,0)[l]{$K_\rho$}}
                        \end{picture}}

\psfrag{Lp}[][]{\begin{picture}(0,0)
            \put(-6,-1){\makebox(0,0)[l]{$L_\rho$}}
                        \end{picture}}

\psfrag{Mp}[][]{\begin{picture}(0,0)
            \put(-16,4){\makebox(0,0)[l]{$M_\rho$}}
                        \end{picture}}

\psfrag{Kbpp}[][]{\begin{picture}(0,0)
            \put(4,-2){\makebox(0,0)[l]{$K^\prime$}}
                        \end{picture}}

\psfrag{Xp}[][]{\begin{picture}(0,0)
            \put(-9,-2){\makebox(0,0)[l]{$X$}}
                        \end{picture}}

\psfrag{Yp}[][]{\begin{picture}(0,0)
            \put(-4,0){\makebox(0,0)[l]{$Y$}}
                        \end{picture}}

\psfrag{Zp}[][]{\begin{picture}(0,0)
            \put(-11,-2){\makebox(0,0)[l]{$Z$}}
                        \end{picture}}

\setlength{\abovecaptionskip}{4pt}%
\setlength{\belowcaptionskip}{-15pt}%

\resizebox{!}{!}{%
  \includegraphics{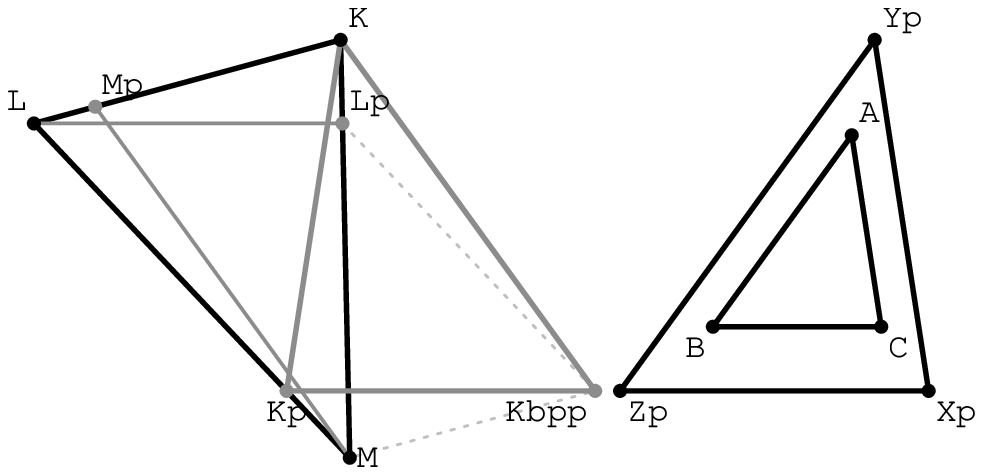}}
    \caption{$T=ABC, V=KLM, \rho =4/5, \CT_{\rho} (V) = XY\!Z$ is directly similar to $T$}
\label{fSBAb}
\end{figure}
\end{example}

\smallskip

In the following corollary we prove that the function $T \mapsto
\gamma_{_T} + i \mu_{_T}$ is a shape function. What we mean by this
is that oriented triangles $T$ and $V$ are directly similar if and
only if $\gamma_{_T} + i \mu_{_T} = \gamma_{_V} + i \mu_{_V}$. To
include equilateral triangles as well, we set $\mu_{_T} = 1$ if and
only if $T$ is equilateral.

\begin{corollary} \label{cshape}
Let $T$ and $V$ be oriented triangles. Then, $T$ and $V$ are directly similar if and only if $\gamma_{_T} = \gamma_{_V}$ and $\mu_{_T} = \mu_{_V}$.
\end{corollary}
\begin{proof}
Let $T$ and $V$ be directly similar. Then, clearly,  $\gamma_{_T} = \gamma_{_V}$. By Corollary~\ref{rho-T-V},  $\CT_{\mu_{_T}}\!(T)$ and $\CT_{\mu_{_T}}\!(V)$ are directly similar wide triangles. Theorem \ref{tit} (\ref{titie}) implies $\mu_{_T} \in \mu_{_V} \Bo\nS$. Since both $\mu_{_T}$ and $\mu_{_V}$ are in $\nI$, they must be equal.

Conversely, assume $\gamma_{_T} = \gamma_{_V}$ and $\mu_{_T} = \mu_{_V}$. Set, $W_T = \CT_{\mu_{_T}}(T)$ and $W_V = \CT_{\mu_{_V}}(V)$. By Theorem~\ref{th1}, $\gamma_{_{W_T}} = \gamma_{_T} = \gamma_{_V} = \gamma_{_{W_V}}$. For a wide isosceles triangle $W$ with sides $a=b<c$, using \eqref{eqcga}, we have
\begin{equation*}
3 (\cos \gamma_{_W})^2 = \frac{\left(2 + \bigl(\frac{c}{b}\bigr)^2\right)^2}{2+\bigl(\frac{c}{b}\bigr)^4}.
\end{equation*}
The last equation has a unique solution for $c/b$ in the interval $(1,2)$. That solution is
\[
\frac{c}{b} =  \sqrt{
 \frac{\sqrt{2}+ 2 \tan \gamma_{_W}}{\sqrt{2} - \tan \gamma_{_W}}}=
 \sqrt{\sqrt{2}\tan\big(\gamma_{_W}+\arctan (1/\sqrt{2})\big)}
 \in (1,2).
\]
Hence, $W_T$ and $W_V$ are directly similar wide triangles.  Since $\mu_{_T} = \mu_{_V}$, Corollary~\ref{rho-T-V} yields that $T$ and $V$ are directly similar.
\end{proof}

\begin{example}
The function $T \mapsto \gamma_{_T} + i \varkappa_{_T}$ is not a shape function since the function $T \mapsto \varkappa_{_T}$ does not have the same values on similar triangles. The function $T \mapsto \gamma_{_T} + i \pI(\varkappa_{_T}\!)$ does take the same values on similar triangles, but it is not a shape function. Indeed, with triangles $T$  and $V$ from Example~\ref{exVsT} we have $\tan \gamma_{_V} = \tan \gamma_{_U} = \sqrt{2}/5$, $\varkappa_{_V} = \varkappa_{_U} = 1/3$ and hence $\pI(\varkappa_{_V}\!) = \pI(\varkappa_{_U}\!) = 1/3$. The last equality explains why $\nM_T = \nM_V$ in Example~\ref{exVsT}. We notice that it can be proved that the restriction of $T \mapsto \gamma_{_T} + i \pI(\varkappa_{_T}\!)$ to the family of increasing (or decreasing) triangles is a shape function.
\end{example}

\begin{remark}
The last corollary implies that two parameters, $\mu \in (0,1/2)$ and  $\gamma \in \bigl(0,\arctan (1/\sqrt{2})\bigr)$, determine an oriented triangle uniquely up to direct similarity. To find sides $a, b, c$ of such a triangle one would need to solve the following system of equations:
\begin{align*}
(\cos \gamma)^2 & = \frac{(a^2 + b^2 + c^2)^2}{3(a^4+b^4+c^4)}, \\
\mu & = \frac{b^2 - a^2}%
 {c^2 - a^2 + \sqrt{a^4 + b^4 + c^4 - a^2 b^2 - b^2 c^2 - c^2a^2}}
\end{align*}
for $a \leq b < c$. Our theory yields a relatively simple family of solutions. They are the sides of the triangle $\CT_\mu (W)$, where $W$ is a wide triangle with the sides
\begin{equation*}
 t = t < t\,F, \quad t>0, \qquad \text{with} \quad F = \sqrt{
 \frac{\sqrt{2}+ 2 \tan \gamma}{\sqrt{2} - \tan \gamma}}.
\end{equation*}
That is,
\begin{align*}
a & = t \sqrt{1-(1-\mu ) \mu F }, \quad
 b   = t \sqrt{\mu F + (1- \mu)^2}, \quad
  c  = t \sqrt{(1-\mu ) F + \mu^2}.
\end{align*}

\end{remark}

\section{Applications and examples}

In the first application of our results we characterize those triangles $T$ whose family of Ceva's triangles contains a right triangle.

\begin{theorem}
Let $T$ be an oriented non-equilateral triangle. There exists $\rho \in \nRc$ such that $\CT_\rho(T)$ is a right triangle if and only if $\tan \omega_{_T} \leq 1/2$.
\end{theorem}
\begin{proof}
Assume that $R=\CT_\rho(T)$ is a right triangle. Then $a^2+b^2=c^2$, and thus, by \eqref{eqcga},
\begin{align*}
(\cos \gamma_{_R})^2 & = \frac{4 c^4}{3((c^2-b^2)^2+b^4+c^4)} = \frac{2}{3\left(1  - \bigl(\tfrac{b}{c}\bigr)^2 + \bigl(\tfrac{b}{c}\bigr)^4\right)} \\
& = \frac{2}{\tfrac{9}{4} + 3 \left(\tfrac{1}{2} - \bigl(\tfrac{b}{c}\bigr)^2\right)^2} \leq \frac{8}{9}.
\end{align*}
Hence, $\tan \gamma_{_R} \geq \sqrt{2}/4$. By \eqref{eqbaca} this is equivalent to $\tan \omega_{_R} \leq 1/2$. Since, by Theorem~\ref{th1}, $\omega_{_T} = \omega_{_R}$, the necessity is proved.

Now, assume that $\tan \omega_{_R} \leq 1/2$, or equivalently, $\tan \gamma_{_T} \geq \sqrt{2}/4$. Let $W = \CT_{\mu_{_T}}(T)$ be the wide isosceles Ceva's triangle in $\CT_{\nR}(T)$. Then, by Theorem~\ref{th1}, $\tan \gamma_{_W} = \tan \gamma_{_T} \geq \sqrt{2}/4$.  If $\tan \gamma_{_T} = \sqrt{2}/4$, then $\cos \gamma_{_T} = 2\sqrt{2}/3$ and $W$ is a right triangle by \eqref{eqcga}.  In this case $\rho = \mu_{_T}$ in the theorem.  Assume now that $\tan \gamma_{_T} > \sqrt{2}/4$. Denote by $F_W$ the ratio of the base to the leg of $W$. Since $\tan \gamma_{_W} = \tan \gamma_{_T}$ using a formula from the proof of  Corollary~\ref{cshape} we have
\[
F_W = \sqrt{
 \frac{\sqrt{2}+ 2 \tan \gamma_{_T}}{\sqrt{2} - \tan \gamma_{_T}}} \in \bigl(\sqrt{2},2\bigr).
\]
Without loss of generality, we can assume that the sides of $W$ are $a = b = 1 < F_W$. We will
now calculate $\tau \in \nI$ such that $\CT_\tau(W)$ is a right triangle. We use the formulas \eqref{Crho}, \eqref{Brho}, \eqref{Arho} and set up the equation for $\tau$: $x_\tau^2 + y_\tau^2 = z_\tau^2$. Solving this equation for $\tau$ gives
\[
\tau = \frac{2-F_W^2+\sqrt{\bigl(F_W^2-2\bigr)\bigl(5 F_W^2-2\bigr)}}{2 F_W^2} \in (0,1/2).
\]
With this $\tau$, $\CT_\tau(W) = \CT_\tau\bigl(\CT_{\mu_{_T}}(T)\bigr)$ is a right triangle. Since, by Corollary~\ref{ctds}, $\CT_\tau\bigl(\CT_{\mu_{_T}}(T)\bigr)$ is directly similar to $\CT_{\varkappa_{_T}\Bo \mu_{_T} \Bo \Ton{\tau} }(T)$, the theorem is proved.
\end{proof}

The second application addresses a slight ambiguity left in Section~\ref{sshape}. All our statements in Section~\ref{sshape} involve the concept of similarity. We are thus interested in finding out whether it is indeed possible to reconstruct a {\it congruent copy} of a given oriented triangle $V$ from an oriented triangle $T$ having the same Brocard angle as $V$. By Theorem~\ref{th1}, one application of the operator $\CT_\rho$ to $T$ with the correct choice of $\rho$ produces a similar copy of $V$, but not necessarily a congruent one. A natural choice here is to iterate sufficiently many times the operators $\CT_\rho$ with possibly different parameters $\rho$. The binary similarity property of the iterations of Ceva's operators is the key to the proof below.

\begin{theorem}\label{th2}
Let $T$ and $V$ be oriented non-equilateral triangles. The triangles $T$ and $V$ have the same Brocard angle if and only if there exist a nonnegative integer $n_0$ and $\rho, \xi \in \nRc$ such that $\CT_\rho\bigl(\CT^{2n_0}_\xi(T)\bigr)$ is directly congruent to $V$.
\end{theorem}
\begin{proof}
Assume that $\omega_{_T} = \omega_{_V}$. By Theorem~\ref{th1} and Remark~\ref{rVsT}, the triangle $V$ is directly similar to $\CT_{\mu_{_T}\BoI\mu_{_V}}\!\!(T)$. Let $\rho = \mu_{_T}\BoI\mu_{_V}$. Then $V$ is directly similar to $\CT_{\rho}(T)$, that is, there exists $l > 0$ such that $V$ and $l \CT_{\rho}(T)$ are directly congruent.

Since $l^{1/n} \to 1$ as $n \to +\infty$, there exists a minimum $n_0 \in \nN$ such that $l^{1/n_0} \geq 3/4$. Also, since $x \mapsto \pn{x}^2$ is a quadratic function with the vertex at $(1/2,3/4)$, there exists a unique $\xi \geq 1/2$ such that $l^{1/n_0} = \pn{\xi}^2$, that is $l  = \pn{\xi}^{2 n_0}$. The last claim in Corollary~\ref{rho-T-V} yields $\CT^{2n_0}_\xi(T) = \pn{\xi}^{2n_0} T = l\, T$. Applying Corollary~\ref{rho-T-V} again we get $\CT_\rho\bigl(\CT^{2n_0}_\xi(T)\bigr) = l\,\CT_\rho(T)$. Since $V$ is directly congruent to $l\,\CT_\rho(T)$, the sufficiency part of our theorem is proved. The necessity part follows by applying Theorem~\ref{th1} $2n_0+1$ times.
\end{proof}

\begin{figure}[ht]

\setlength{\abovecaptionskip}{10pt}%
\setlength{\belowcaptionskip}{4pt}%

\psfrag{Z}[][]{\begin{picture}(0,0)
            \put(-3,-14){\makebox(0,0)[l]{$Z$}}
                        \end{picture}}

\psfrag{Y}[][]{\begin{picture}(0,0)
            \put(-3,0){\makebox(0,0)[l]{$Y$}}
                        \end{picture}}

\psfrag{X}[][]{\begin{picture}(0,0)
            \put(-9,-2){\makebox(0,0)[l]{$X$}}
                        \end{picture}}

\psfrag{B}[][]{\begin{picture}(0,0)
            \put(0,-2){\makebox(0,0)[l]{$B$}}
                        \end{picture}}

\psfrag{Ao}[][]{\begin{picture}(0,0)
            \put(-6,-3){\makebox(0,0)[l]{$A_\varkappa$}}
                        \end{picture}}

\psfrag{A}[][]{\begin{picture}(0,0)
            \put(-10,2){\makebox(0,0)[l]{$A$}}
                        \end{picture}}

\psfrag{C}[][]{\begin{picture}(0,0)
            \put(-9,-2){\makebox(0,0)[l]{$C$}}
                        \end{picture}}

\psfrag{Aop}[][]{\begin{picture}(0,0)
            \put(-10,-12){\makebox(0,0)[l]{$A_\varkappa'$}}
                        \end{picture}}

\psfrag{D}[][]{\begin{picture}(0,0)
            \put(-8,2){\makebox(0,0)[l]{$D$}}
                        \end{picture}}

\resizebox{!}{!}{%
  \includegraphics{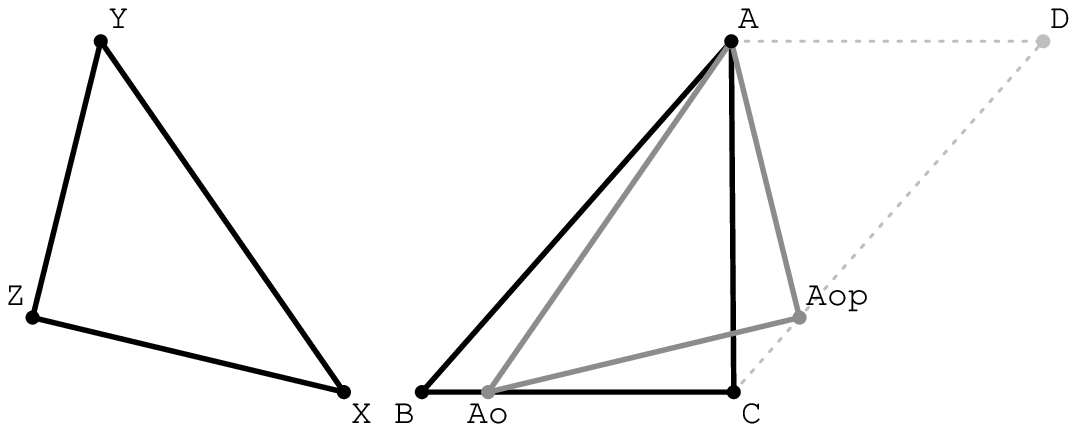}}
    \caption{$\CT_{\varkappa} (T) = XY\!Z$ is directly similar to $T=ABC$}
\label{f2}

\psfrag{Zm}[][]{\begin{picture}(0,0)
            \put(-2,-14){\makebox(0,0)[l]{$Z$}}
                        \end{picture}}

\psfrag{Zn}[][]{\begin{picture}(0,0)
            \put(-6,-10){\makebox(0,0)[l]{$M$}}
                        \end{picture}}

\psfrag{Ym}[][]{\begin{picture}(0,0)
            \put(-8,0){\makebox(0,0)[l]{$Y$}}
                        \end{picture}}

\psfrag{Yn}[][]{\begin{picture}(0,0)
            \put(-7,0){\makebox(0,0)[l]{$L$}}
                        \end{picture}}

\psfrag{Xm}[][]{\begin{picture}(0,0)
            \put(-32,-3){\makebox(0,0)[l]{$X = K$}}
                        \end{picture}}

\psfrag{Am}[][]{\begin{picture}(0,0)
            \put(-6,-4.5){\makebox(0,0)[l]{$A_\mu$}}
                        \end{picture}}

\psfrag{Amp}[][]{\begin{picture}(0,0)
            \put(-10,-12){\makebox(0,0)[l]{$A_\mu'$}}
                        \end{picture}}

\psfrag{An}[][]{\begin{picture}(0,0)
            \put(-6,-4){\makebox(0,0)[l]{$A_\nu$}}
                        \end{picture}}

\psfrag{Anp}[][]{\begin{picture}(0,0)
            \put(-10,-12){\makebox(0,0)[l]{$A_\nu'$}}
                        \end{picture}}

\psfrag{B}[][]{\begin{picture}(0,0)
            \put(-2,-3){\makebox(0,0)[l]{$B$}}
                        \end{picture}}

\psfrag{C}[][]{\begin{picture}(0,0)
            \put(-9,-3){\makebox(0,0)[l]{$C$}}
                        \end{picture}}

\resizebox{!}{!}{%
  \includegraphics{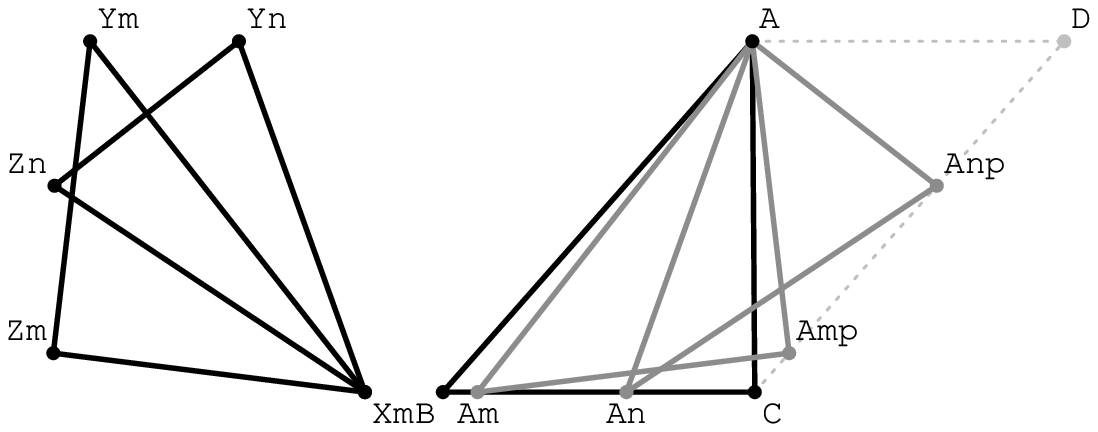}}
    \caption{The isosceles triangles $\CT_\mu (T) = XY\!Z$ and $\CT_\nu (T) = KLM$}
\label{f3}

\psfrag{Zp}[][]{\begin{picture}(0,0)
            \put(-3,-12){\makebox(0,0)[l]{$Z$}}
                        \end{picture}}

\psfrag{Zq}[][]{\begin{picture}(0,0)
            \put(-7,-9){\makebox(0,0)[l]{$M$}}
                        \end{picture}}

\psfrag{Yp}[][]{\begin{picture}(0,0)
            \put(-7,0){\makebox(0,0)[l]{$Y$}}
                        \end{picture}}

\psfrag{Yq}[][]{\begin{picture}(0,0)
            \put(-7,0){\makebox(0,0)[l]{$L$}}
                        \end{picture}}

\psfrag{Xp}[][]{\begin{picture}(0,0)
            \put(-32,-1){\makebox(0,0)[l]{$X = K$}}
                        \end{picture}}

\psfrag{Ap}[][]{\begin{picture}(0,0)
            \put(-6,-5){\makebox(0,0)[l]{$A_\xi$}}
                        \end{picture}}

\psfrag{App}[][]{\begin{picture}(0,0)
            \put(-10,-12){\makebox(0,0)[l]{$A_\xi'$}}
                        \end{picture}}

\psfrag{Aq}[][]{\begin{picture}(0,0)
            \put(-6,-5){\makebox(0,0)[l]{$A_{\imath(\xi)}$}}
                        \end{picture}}

\psfrag{Aqp}[][]{\begin{picture}(0,0)
            \put(-10,-12){\makebox(0,0)[l]{$A_{\imath(\xi)}'$}}
                        \end{picture}}

\psfrag{B}[][]{\begin{picture}(0,0)
            \put(0,-3){\makebox(0,0)[l]{$B$}}
                        \end{picture}}

\psfrag{C}[][]{\begin{picture}(0,0)
            \put(-9,-3){\makebox(0,0)[l]{$C$}}
                        \end{picture}}

\setlength{\abovecaptionskip}{10pt}%
\setlength{\belowcaptionskip}{-15pt}%

\resizebox{!}{!}{%
  \includegraphics{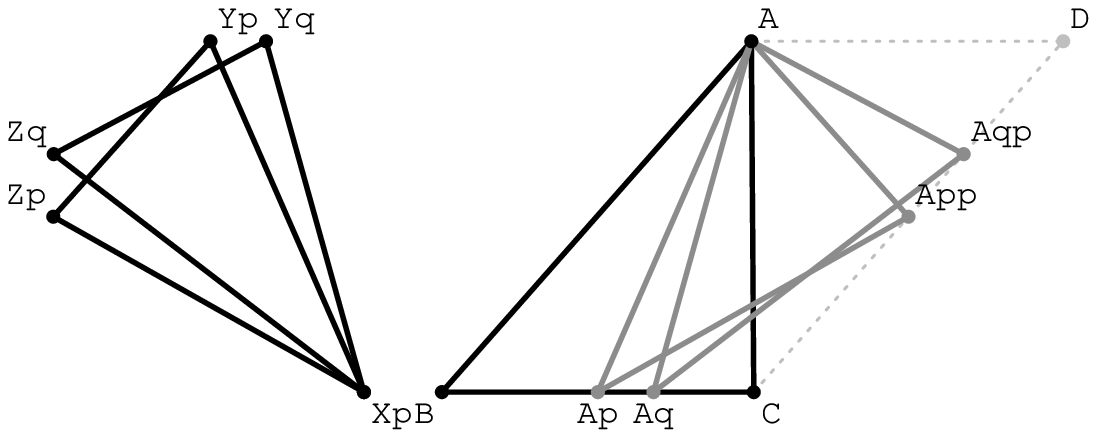}}
    \caption{$\CT_\xi (T)=XY\!Z$ is reversely similar to $\CT_{\imath(\xi)} (T)=KLM$}
\label{f4}
\end{figure}

\begin{example}
We illustrate Theorem~\ref{th2} with a scaled median triangle. Let  $V = (1/4)\CT_{1/2}(T)$, that is the triangle that is similar to the median triangle with the ratio $1/4$. By Corollary~\ref{comega}, the only triangles in $\CT_\nR(T)$ similar to  $\CT_{1/2}(T)$ are $\CT_{-1}(T)$ and $\CT_2(T)$. The ratio of similarity of each of them to $\CT_{1/2}(T)$ is $2$. Hence, $V$ is not a Ceva's triangle of~$T$. But we can apply  Theorem~\ref{th2}, with $l = 1/4$ and $\rho = 1/2$ introduced in its proof.  Since $(3/4)^4 > 1/4$ and $(3/4)^5 < 1/4$, we have $n_0 = 5$. We choose $\xi$ such that  $1-\xi+\xi^2 = \sqrt[5]{1/4}$, that is
\[
\xi = \frac{1}{2}\left(1-\sqrt{2\,2^{3/5}-3}\right)\approx 0.4114 \quad  \text{or} \quad \xi = \frac{1}{2}\left(1+\sqrt{2\,2^{3/5}-3}\right) \approx 0.5886.
\]
Then $\pn{\xi}^2 = 2^{-2/5}$ and, by Corollary~\ref{rho-T-V},
\[
\CT_\xi^{10}(T) = \pn{\xi}^{10}\,T = \bigl(2^{-2/5}\bigr)^5 \, T = \frac{1}{4}\, T.
\]
By Corollary~\ref{rho-T-V} again,
\[
\CT_{1/2}\bigl(\CT_\xi^{10}(T)\bigr) = \frac{1}{4} \, \CT_{1/2}(T).
\]
Thus, a composition of eleven Ceva's operators reconstructs  $\frac{1}{4} \CT_{1/2}(T)$ from~$T$.
\end{example}

\begin{example}
We conclude the paper by illustrating the fundamental numbers
$\varkappa_{_T}, \mu_{_T}, \nu_{_T}$ associated with a given
oriented triangle $T$, Theorems~\ref{CJT} and \ref{th2} and Remark~\ref{rVsT}. Let $T = ABC$ be the increasing triangle with sides $a = 8$, $b = 9$ and $c = 12$.  We calculate $\varkappa = \varkappa_{_T}=17/80$; see Figure~\ref{f2}, where $\CT_\varkappa (T) = XY\!Z$ is directly similar to $T$.

Furthermore, we have $\mu = \mu_{_T}=1/9$ and $\nu = \nu_{_T} =
10/17$. In agreement with Theorem~\ref{tit}, $\CT_\mu (T)= XY\!Z$ is a
wide triangle, $\CT_\nu (T)=KLM$ is a narrow triangle and $\mu <
\varkappa <\nu$, see Figure~\ref{f3}. We also observe that the bases
$XY$ and $LM$ of the isosceles triangles are perpendicular. It can
be confirmed that this is always true by computing the dot product
of the vectors $\overrightarrow{AA}_\mu = \overrightarrow{AB} + \mu
\overrightarrow{BC}$ and $\overrightarrow{CC}_\nu =
\overrightarrow{C\!A} + \nu \overrightarrow{AB}$ in terms of the sides
$a,b,c$ and using the formulas for $\mu$ and $\nu$ given in
Remark~\ref{rmuf}.

The median triangle of $T$, that is Ceva's triangle $\CT_{1/2}(T)$ is calculated using \eqref{Crho}, \eqref{Brho}, \eqref{Arho} with $\rho = 1/2\in \nM_{_T}$, to be
\[
\CT_{1/2} (T) =
\left(\sqrt{\frac{73}{2}},\frac{\sqrt{335}}{2},\sqrt{\frac{193}{2}}\right).
\]
With $\imath$ defined in \eqref{involution} we get $\imath(1/2)= 97/143\in \nI\setminus \nM_{_T}$ and calculate
\[
\CT_{97/143} (T)= \left( \frac{73\sqrt{146}}{143},
   \frac{73\sqrt{386}}{143},
   \frac{73\sqrt{335}}{143}\right).
\]
By Theorem~\ref{CJT} triangles $\CT_{1/2} (T)$ and $\CT_{97/143} (T)$ are reversely similar with the ratio of similarity  $143/146$. This is illustrated in Figure~\ref{f4}.
 \end{example}

\end{document}